\documentclass[preprint,12pt]{elsarticle}
\usepackage{doi}
\usepackage{lineno,hyperref}
\usepackage{amssymb,amsfonts}
\usepackage{amsmath,amsthm}
\usepackage{algorithm,algorithmic}
\usepackage{caption}
\usepackage{graphicx,subfigure}
\usepackage{float,rotating}
\usepackage{booktabs,array}
\usepackage{multirow,multicol}
\usepackage{lscape}
\usepackage{epstopdf}
\usepackage{color}
\usepackage{setspace}
\usepackage{natbib}
\usepackage{bm}
\usepackage{url,doi}

\newtheorem{lemma}{Lemma}

\usepackage{amssymb}
\usepackage{amsmath}

\begin{document}

\begin{frontmatter}



\title{Splitting-based randomized dynamical low-rank approximations for stiff matrix differential equations}

\author[address1]{Zi Wu}
\ead{wuziiii@163.com}

\author[address1]{Yong-Liang Zhao\corref{correspondingauthor}}
\cortext[correspondingauthor]{Corresponding author}
\ead{ylzhaofde@sina.com}

\author[address2]{Xian-Ming Gu}
\ead{guxianming@live.cn, guxm@swufe.edu.cn}


\affiliation[address1]{organization={School of Mathematical Sciences, Sichuan Normal University},
	city={Chengdu},
	postcode={610068},
	state={Sichuan},
	country={P.R. China}}

\affiliation[address2]{organization={School of Mathematics, Southwestern University of Finance and Economics},
	city={Chengdu},
	postcode={611130},
	state={Sichuan},
	country={P.R. China}}

%

\begin{abstract}
In the fields of control theory and machine learning, dynamic low-rank approximations for time-dependent matrices have received substantial attention. This work considers a low-rank solution of the large-scale semilinear stiff matrix differential equation from such fields. We first split such the equation into a stiff linear subproblem and a nonstiff nonlinear subproblem. Then, a low-rank exponential integrator is applied to the linear subproblem. Two randomized low-rank approaches are developed for the nonlinear subproblem. Furthermore, we extend these low-rank methods to rank-adaptation scenarios. Through rigorous validation on canonical stiff matrix differential problems, including spatially discretized Allen-Cahn equations and differential Riccati equations, numerical experiments demonstrate that the proposed methods attain the expected convergence orders and exhibit strong robustness and accuracy.
\end{abstract}

%


\end{frontmatter}



\section{Introduction}
\label{sec1}

Low-rank approximation of large-scale matrices has demonstrated significant value in various fields, including bioinformatics \cite{kapur2016gene}, machine learning \cite{sainath2013low,zhang2021general}, mathematical finance \cite{stoica2023low}, control theory \cite{benner2018,lancaster1995}, image compression \cite{andrews1976singular}, and natural language processing \cite{deerwester1990indexing}. In this work, we seek low-rank approximate solutions to large-scale semilinear stiff matrix differential equations of the following form:
\begin{equation}\label{eq1.1}
	\dot{X}(t)=AX(t)+X(t)A^{\top}+F(t,X(t)),
	\quad X(t_{0})=X_{0},
\end{equation}
where $X(t)\in\mathbb{R}^{m\times n},F:[t_{0},\infty)\times\mathbb{R}^{m\times n}\to\mathbb{R}^{m\times n}$ is nonlinear, and $X_0$,   $A\in\mathbb{R}^{m\times n}$ are time-independent.
Such equations arise in a variety of large-scale applications. For instance, spatial semi-discretizations of nonlinear PDEs, such as the Allen--Cahn equation, usually lead to stiff linear terms generated by discrete differential operators. Differential Riccati equations from control theory provide another important example, where large-scale linear dynamics are coupled with nonlinear quadratic matrix terms. These problems motivate the development of efficient low-rank methods that can handle stiffness, nonlinearity, and high dimensionality simultaneously.



In recent years, efficient numerical solutions for large-scale stiff matrix differential equations have garnered significant attention in the fields of physics, chemistry, and engineering \cite{jang2023sparse,weng2022multiscale}. The primary challenge lies in balancing computational efficiency with robustness against stiffness and high-dimensional nonlinearity. The traditional Lie-Trotter splitting method addresses the linear part analytically by separating the stiff and nonstiff terms, with recent extensions demonstrating efficacy for nonlinear fractional PDEs \cite{zhao2021low}. Then an exponential integrator is employed   \cite{hairer1993sp,ostermann2019convergence}. However, a dynamic low-rank approximation (DLRA) of the nonlinear terms still depends on the projective-splitting integrator \cite{lubich2014projector}. It faces challenges related to sensitivity to small singular values and bottlenecks in computational efficiency \cite{kieri2016discretized}. Although the DLRA significantly reduces dimensionality through manifold constraints \cite{koch2007dynamical}, its classical implementation is constrained by explicit step size limitations and fixed rank assumptions \cite{kieri2019projection}. Subsequent improvements, such as the adaptive BUG method \cite{carrel2025randomized,ceruti2022rank}, enhance the stability but have limited scalability for high-dimensional nonlinear problems. The emergence of randomized numerical linear algebra introduces new ideas for low-rank decomposition, including the randomized singular value decomposition (RSVD) proposed by Halko et al. \cite{halko2011finding} and the generalized Nyström method developed by Tropp et al. \cite{tropp2017practical}. These methods accelerate matrix approximation through sketching techniques. Additionally,  Carrel \cite{carrel2025randomized} extends these approaches to dynamic scenarios by developing the dynamic RSVD (DRSVD) and dynamic generalized Nyström (DGN) methods. These methods demonstrate low variance and high parallelism in time-dependent systems \cite{lam2024randomized}. However, most existing studies concentrate on static or weakly stiff scenarios, and the deep integration of the splitting framework with randomized dynamic low-rank (RDLR) methods has not been thoroughly investigated.

In this work, inspired by \cite{carrel2025randomized}, to address the stiffness, nonlinearity, and unbearably large size of Eq.~\eqref{eq1.1}, splitting-based RDLR approximations are designed. Our strategy not only inherits the robustness of the splitting method in handling stiffness \cite{sutti2024implicit} but also enhances the computational efficiency of the nonlinear terms through randomization techniques. To balance the accuracy and computational efficiency of the proposed methods, especially when the optimal rank is unknown beforehand and varies greatly over time, we extend them to rank-adaptive settings. It also provides a novel solving paradigm for higher-dimensional stiff systems, which has both high efficiency and stability. The main contributions can be summarized as follows:
\begin{itemize}
	\item[(a)] A novel RDLR method based on an operator splitting strategy has been developed, specifically tailored to the structural characteristics of Eq.~\eqref{eq1.1}. This effectively decouples and handles the linear stiff term and the nonlinear term, enabling efficient computation with longer time and steps simulations than \cite{carrel2025randomized}. For instance, the Allen-Cahn equation inherently favors these features.
	\item[(b)] The linear stiff term is integrated exactly using matrix exponentials \cite{hochbruck2010exponential}. This overcomes explicit step size restrictions inherent in standard methods sizes \cite{kassam2005fourth,leveque2007finite}, ensuring unconditional stability.
	\item[(c)] The nonlinear term is resolved efficiently via the RDLR method \cite{carrel2025randomized}. This approach bypasses the computational bottlenecks associated with high-dimensional matrix operations.
	\item[(d)] The proposed framework unifies computational efficiency, numerical accuracy, and stiffness stability. It is particularly well-suited for large-scale scientific computations, including high-dimensional PDEs, control theory \cite{benner2018,lancaster1995} (e.g., Riccati equations), and dynamical systems \cite{koch2007dynamical}.
\end{itemize}

The rest of this paper is organized as follows. Section \ref{sec2} provides details of our proposed methods. We first establish the Lie-Trotter and Strang splitting formats for Eq.~\eqref{eq1.1}, respectively. Then, based on these splittings, splitting-based RDLR approximations are designed. Section \ref{sec3} extends these approaches to the rank-adaptive scenario. Numerical validation in Section \ref{sec4} includes convergence tests and large-scale simulations, demonstrating the excellent numerical performance. Section \ref{sec5} presents the main conclusions and outlines directions for future research. The pseudocode of some key algorithms discussed in this paper is provided in Appendix.


\section{Splitting-based randomized dynamical low-rank approximations}
\label{sec2}

By utilizing the variation-of-constants formula \cite{collegari2018linear}, the exact solution of Eq.~\eqref{eq1.1} can be expressed as follows:
\begin{equation*}
	X(t)=\mathrm{e}^{(t-t_0)A}X(t_0)\mathrm{e}^{(t-t_0)A^{\top}}+\int_{t_0}^t\mathrm{e}^{(t-s)A}F(s,X(s))\mathrm{e}^{(t-s)A^{\top}}\mathrm{d}s,
\end{equation*}
for $t_{0}\leq t\leq T$.

The stiffness of Eq.~\eqref{eq1.1} limits the applicability of most standard numerical integrators \cite{hochbruck2010exponential}. Explicit methods must satisfy a CFL-like condition and require the use of tiny time step sizes \cite{darwish2016finite}, while implicit methods result in significant computational costs. Consequently, we propose RDLR approximations based on the splitting method.

\subsection{Splitting into two subproblems}
\label{sec2.1}

Taking the structure of \eqref{eq1.1} into account, the key idea is the use of splitting method.  For foundational theory and modern developments in operator splitting, see \cite{hochbruck2010exponential,mclachlan2002splitting}.

Now, we split \eqref{eq1.1} into two terms: the stiff (linear) part and the nonstiff (nonlinear) part. This yields the following two subproblems:

\begin{equation}\label{eq2.1}
	\dot{M}(t)=AM(t)+M(t)A^{\top},\quad\quad M(t_{0})=M_{0},
\end{equation}
and
\begin{equation}\label{eq2.2}
	\dot{N}(t)=F(t,N(t)),\quad\quad N(t_0)=N_0,
\end{equation}
where $t\in[t_{0},T]$. The principal advantage of this splitting strategy stems from its ability to precisely decouple stiffness dynamics: the stiff component is entirely encapsulated within the first subproblem \eqref{eq2.1}, whereas the second subproblem  \eqref{eq2.2} is a nonstiff ODE. This separation achieves dual computational benefits: (i) it expands the admissible choices of numerical integrators for the nonstiff subsystem, and (ii) enables systematic derivation of splitting schemes with arbitrary convergence orders by reconfiguring the operator composition in \eqref{eq1.1}.

For a given positive integer $\hat{M}$, the time step size $\tau$ is defined as $\tau=\frac{T-t_0}{\hat{M}}$. Then, the time interval $[t_{0},T]$ is decomposed as $[t_0,T]=\bigcup\limits_{k=1}^{\hat{M}} [t_0, t_0 + k\tau]$. We denote the solutions of the subproblems (2) and (3) with initial values $M_0$ and $V_0$ by $\Phi_\tau^A(M_0)$ and $\Phi_\tau^F(N_0)$, respectively. 
The simplest splitting method for solving Eq.~\eqref{eq1.1} is the Lie-Trotter splitting scheme with step size $\tau$:
\begin{equation*}
	\mathcal{L}_\tau:=\Phi_\tau^A\circ\Phi_\tau^F.
\end{equation*}
With this notation, we can obtain an approximation $X_{\tau}^{L}$ of the solution to $X(t)$ \eqref{eq1.1} at $t = t_0+\tau$, i.e.,
\begin{equation*}
	X_{\tau}^{L}=\mathcal{L}_\tau(X_{0})=\Phi_\tau^A\circ\Phi_\tau^F(X_{0}),
\end{equation*}
where $N_{0} = X_{0}$. This is a first-order method, and by the symmetric formula, we can have a second-order method called Strang splitting: 
\begin{equation*}
	S_{\tau}:=\Phi_{\tau/2}^{A}\circ\Phi_{\tau}^{F}\circ\Phi_{\tau/2}^{A}.
\end{equation*}
Similarly, the solution of \eqref{eq1.1} at $t_0 + \tau$ is approached as: 
\begin{equation*}
	X_{\tau}^{S}=S_\tau(X_{0})=\Phi_{\tau/2}^A\circ\Phi_\tau^F\circ\Phi_{\tau/2}^A(X_{0}),
\end{equation*}
where $M_{0} = X_{0}$.

Based on these two splitting methods, a low-rank solution of \eqref{eq1.1} is investigated in the next. 


\subsection{The dynamical low-rank integration}
\label{sec2.2}

Denoting the manifold of rank-$r$ matrices by 
\begin{equation*}
	\mathcal{M}_r=\{M\in\mathbb{R}^{m\times n}\mid\text{rank}(M)=r\}.
\end{equation*}
Our aim is to find a low-rank approximation \(Y \in \mathcal{M}_r\) $(r\ll\min\{m,n\})$ to the solution of \eqref {eq1.1}.
Under the splitting strategy introduced in Section \ref{sec2.1}, our primary objective turns to seeking low-rank approximations \(W \in \mathcal{M}_r\) and \(Z \in \mathcal{M}_r\) for subproblems \eqref {eq2.1} and \eqref {eq2.2}, respectively.

We first concentrate on the stiff subproblem \eqref {eq2.1}. For any \(Y \in \mathcal{M}_r\), we observe that \(AY + YA^{\top} \in \mathcal{T}_Y\mathcal{M}_r\), where $\mathcal{T}_Y\mathcal{M}_r$ is the tangent space of $\mathcal{M}_r$ at the current state $Y$ \cite{helmke2012optimization}. This implies that if the initial condition lies in \(\mathcal{M}_r\), the solution of \eqref {eq2.1} remains within \(\mathcal{M}_r\) for all times. In other words, subproblem \eqref {eq2.1} preserves the rank of its solution. 

The exact solution of this homogeneous problem at $t_0+\tau$ within initial value $M(t_{0})=M_{0}$ is given by 
\begin{equation*}
	M(t_0+\tau)=\mathrm{e}^{\tau A}M_{0}\mathrm{e}^{\tau A^{\top}}.
\end{equation*}
If $M_{0}$ is given in rank-$r$ matrix factorization $USV^{\top}$, where $U\in\mathbb{R}^{m\times r}$, $V\in\mathbb{R}^{n\times r}$ have orthonormal columns, $S\in\mathbb{R}^{r\times r}(r\ll\min\{m,n\})$ is a nonsingular matrix. Subsequently, we obtain a rank-$r$ approximation $W_{\tau}$ to the exact solution $M(t_0+\tau)=\Phi_\tau^A(M_{0})$ by
\begin{equation*}
	W_{\tau}=\mathrm{e}^{\tau A}USV(\mathrm{e}^{\tau A}V)^{\top}.
\end{equation*}
For the matrix exponential-matrix multiplications $\mathrm{e}^{\tau A}U$ and $\mathrm{e}^{\tau A}V$, several high performance algorithms exist, such as Taylor series expansion \cite{al2011computing}, Leja point-based interpolation \cite{Caliari2016}, and Krylov subspace approximations \cite{guttel2013rational,saad1992analysis}. In order to maintain orthogonality, the matrices $\mathrm{e}^{\tau A}U$ and $\mathrm{e}^{\tau A}V$ can undergo the QR decomposition to enforce orthogonality. Specifically, we factorize $\mathrm{e}^{\tau A}U=\widetilde{U}R_1$ where $\widetilde{U}\in\mathbb{R}^{m\times r}$ has orthonormal columns and $R_1\in\mathbb{R}^{r\times r}$ is an upper triangular matrix. Analogously, we decompose $\mathrm{e}^{\tau A}V=\widetilde{V}R_2$ yields another orthonormal basis $\widetilde{V}\in\mathbb{R}^{n\times r}$ with upper triangular matrix $R_2\in\mathbb{R}^{r\times r}$. To preserve the rank-$r$ matrix structure, the core matrix becomes $\widetilde{S}=R_1SR_2^\top$, resulting in the low-rank approximation $W_\tau=\widetilde{U}\widetilde{S}\widetilde{V}^\top$ of the solution to \eqref{eq2.1}. 

For the nonlinear subproblem \eqref{eq2.2}, we employ the RDLR methods proposed in \cite{carrel2025randomized}. These new methods primarily involve two steps: a postprocessing step and a range estimation step (see details in Section \ref{sec2.3}). Firstly, using a dynamical rangefinder method to get an orthonormal matrix $Q\in\mathbb{R}^{m\times r}$ ($r\leq\min(m,n)$) such that
\begin{equation*}
	\left\| N(t) - QQ^{\top} N(t) \right\|_F \stackrel{!}{=} \min.
\end{equation*}
Then, we utilize the DRSVD (or DGN) method to obtain a low-rank approximation $Z \in \mathcal{M}_r$ to \eqref{eq2.2}. The author \cite{carrel2025randomized} has demonstrated that the new method achieves high accuracy, low variance, and significant robustness. After applying the RDLR to \eqref{eq2.2}, a low-rank approximation of \(N(t)\) at \(t = t_0+\tau\) with a low-rank initial value $N_0$ is 
\begin{equation*}
	Z_\tau = \tilde{\Phi}_\tau^F(N_0).
\end{equation*}
Here, \(\tilde{\Phi}_\tau^F\) represents the low-rank solution for \eqref{eq2.2}. 

By combining the low-rank solution flow \(\tilde{\Phi}_\tau^F\) of subproblem  \eqref{eq2.2} and the exact flow \(\Phi_\tau^A\) of subproblem  \eqref{eq2.1}, we can get a rank-$r$ approximate solution $X^L_1$ of Eq.~\eqref{eq1.1} at $t = t_0+\tau$. This is denoted as
\begin{equation}\label{eq:star1}
	\tilde{\mathcal{L}_\tau}:=\Phi_\tau^A\circ\tilde{\Phi}_\tau^F,
\end{equation}
for low-rank Lie-Trotter splitting.
Starting with \(Y_0\), a rank-$r$ approximation of the initial value $X_0$. The rank-\(r\) approximation of \eqref{eq1.1} at time \(t_0+\tau\) is written as
\[
X^\mathcal{L}_1=\tilde{\mathcal{L}_\tau}(Y_0)=(\Phi_\tau^A\circ\tilde{\Phi}_\tau^F)(Y_0).
\]

Similarly, we denote as
\begin{equation*}
	\tilde{\mathcal{S}_\tau}:=\Phi_{\tau/2}^A\circ\tilde{\Phi}_\tau^F\circ\Phi_{\tau/2}^A,
\end{equation*}
for low-rank Strang splitting. Given a rank-$r$ approximation $Y_0$ of the initial value $X_0$, the low-rank approximation of \eqref{eq1.1} at $t_{0}+\tau$ is:
\begin{equation*}
	X^\mathcal{S}_1=\tilde{\mathcal{S}_\tau}(Y_0)=\left(\Phi_{\tau/2}^A\circ\tilde{\Phi}_\tau^F\circ\Phi_{\tau/2}^A\right)(Y_0).
\end{equation*}
To ensure the comprehensiveness and readability of our work, in the next section, we will elaborate how to find a low-rank solution of the nonlinear subproblem \eqref{eq2.2}.


\subsection{Two randomized dynamical low-rank methods for Eq.~\eqref{eq2.2}}\label{sec2.3}
The RDLR framework proposed in \cite{carrel2025randomized} addresses large-scale matrix differential equations through novel adaptations of randomized numerical linear algebra techniques. Its framework is: firstly, a dynamical randomized rangefinder estimates the evolving range subspace; subsequently, two distinct post-processing strategies yield specialized time integrators, i.e., DRSVD and DGN. The RDLR method exhibits distinctive advantages over conventional dynamical low-rank approximation approaches \cite{martinsson2020randomized,schafer2021efficient}: they behave without step size restrictions, and demonstrate inherent robustness. Before propose our two RDLR methods for the nonlinear subproblem \eqref{eq2.2}, we firstly establish a dynamical randomized rangefinder.


To investigate the initial response characteristics of \eqref{eq2.2}, we focuses on analyzing the physical processes within the first time step $[t_0,t_0+\tau]$. The aim of the dynamical rangefinder is to compute an orthonormal basis $Q_\tau\in \mathbb{R}^{m\times(r+p)}$, where $p\geq0$ is the oversampling parameter, such that the column space of $Q_\tau$ approximates the range of $N(t_0+\tau)$. Equivalently, the projection error
\begin{equation*}
	\left\lVert N(t_0+\tau) - Q_\tau Q_\tau^{\top} N(t_0+\tau) \right\rVert_F
\end{equation*}
is expected to be small. Here $||\cdot||_F$ is the Frobenius norm.

A solution to the above optimization problem can be obtained via the following two steps:

1. Let $\Omega \in \mathbb{R}^{n\times(r+p)}$ be a Gaussian random matrix. Instead of explicitly forming $N(t_0+\tau)\Omega$, we evolve the sketched variable $B(t)\in\mathbb{R}^{m\times(r+p)}$ by solving the differential equation:
\begin{equation*}
	\dot{B}(t)
	=
	F\left(t, B(t)(\Omega^\top\Omega)^{-1}\Omega^\top\right)\Omega,
	\qquad
	B(t_0)=N_0\Omega .
\end{equation*}

2. After obtaining $B(t_0+\tau)$, extract the orthonormal basis by
\begin{equation*}
	Q_\tau = \operatorname{orth}\bigl(B(t_0+\tau)\bigr).
\end{equation*}
This can be done, for example, by the QR decomposition.

When higher accuracy is required, the estimated subspace can be further improved by performing $q\geq0$ power iterations through alternating range and co-range updates. The detailed derivation and analysis can be found in \cite{carrel2025randomized}, and the complete procedure is summarized in Algorithm \ref{alg:dyn_rangefinder}.  Now, we are in position to propose two RDLR methods for solving \eqref{eq2.2}.

\subsubsection{The dynamical randomized SVD}\label{sec2.3.2} 
A rank-\(r\) approximation of the initial data \(N_0\) is represented as $N_0 \approx U_0S_0V_0^\top$, where $U_0$ and $V_0$ are two column orthogonal matrices, and
\(S_0\in\mathbb{R}^{r\times r}\) is nonsingular. When this factorisation is computed by a truncated SVD, \(S_0\) is diagonal with positive entries.
The dynamic rangefinder identifies an orthogonal matrix $Q_\tau$ such that $N(t_0 + \tau) \approx Q_\tau Q_\tau^{\top} N(t_0 + \tau)$. This step effectively captures the structure of the low-rank subspaces of $N(t)$ within the first time interval $[t_0,t_0+\tau]$. In order to ensure that the matrix $Q$ includes foundational information along with dynamically evolving orthogonal bases, we combine the initial basis \(U_0\) with the dynamic basis \(Q_\tau\) to construct the extended orthogonal basis \(Q\), defined as 
\begin{equation*}
	Q = \text{orth}([U_0, Q_\tau]).
\end{equation*}
This leads to a novel approximation $N(t) \approx Q Q^{\top} N(t)$ $(t\in[t_0,t_0+\tau])$. It reduces the computational complexity by projecting the original high-dimensional nonlinear equations into the lower-dimensional subspace defined by the matrix $Q$. The projection coefficient matrix $C(t)$ is defined below to derive $n\times(2r+p)$ the differential equation 
\begin{equation}\label{eq2.7}
	\dot{C}(t) = F(t,QC(t)^{\top})^{\top}Q, \quad C(t_0) = N_0^{\top}Q.
\end{equation}

If we apply a SVD to \(C(t_0 + \tau)^{\top}\): 
\begin{equation*}
	\tilde{U}_\tau\Sigma_\tau V_\tau^{\top}=\text{svd}(C(t_0+\tau)^{\top}).
\end{equation*}
Then, we have \(N(t_0 + \tau) \approx N_1 = Q \tilde{U}_\tau \Sigma_\tau V_\tau^{\top}\). Finally, the resulting low-rank decomposition is truncated to the target rank $N_1=\mathcal{T}_r(Q\widetilde{U}_\tau\Sigma_\tau V_\tau^{\top})$. 

The computational complexity can be significantly reduced through the use of low-rank decomposition techniques. In particular, as established in \cite[Lemma 3.1]{carrel2025randomized}, when the system satisfies the regularity condition $\text{Range}(N(t))\subset\text{Range}([U_0,Q_\tau])=\text{Range}(Q)$, for all $t\in[t_0,t_0+\tau]$, Eq.~\eqref{eq2.7} becomes exact reconstruction of the solution.

The DRSVD combines stochastic projection with dynamic evolution by estimating the column space through a dynamic rangefinder, enables efficient low-rank approximation of time-dependent matrix differential equations. The detailed procedure is summarized in Algorithm \ref{alg:drsvd}.

\subsubsection{The dynamical generalised Nystr\"{o}m method}\label{sec2.3.3}

The DGN method provides an alternative approach for solving the nonlinear subproblem~\eqref{eq2.2}. 
In contrast to DRSVD, which mainly relies on an estimated column space, 
DGN constructs a low-rank approximation by dynamically estimating both 
the column space and the row space of $N(t)$ through a dynamical rangefinder and a dynamical co-rangefinder.

Similar to Section \ref{sec2.3.2}, for the rank-r approximation of the initial value $N_0 \approx U_0S_0V_0^\top$, the DGN method proceeds as follows with oversampling parameters $p,\ell \geq 0$, and power iterations $q \geq 0$. 
It worth mentioning that, in here, $S_0$ is not required to be diagonal; if the initial factorisation is obtained from a truncated SVD, then $S_0$ is diagonal with positive singular values.


1. Compute the dynamical rangefinder basis 
$\tilde{Q}_1$ for the column space and the dynamical co-rangefinder basis 
$\tilde{Q}_2$ for the row space:
\begin{align*}
	\tilde{Q}_1 
	&= \operatorname{DynamicalRangefinder}(N_0, \tau, F, r, p, q), \\
	\tilde{Q}_2 
	&= \operatorname{DynamicalCoRangefinder}(N_0, \tau, F, r, p+\ell, q).
\end{align*}

2. Augment the initial bases and form the 
orthonormal matrices:
\begin{equation*}
	Q_1 = \operatorname{orth}([U_0, \tilde{Q}_1]), 
	\qquad 
	Q_2 = \operatorname{orth}([V_0, \tilde{Q}_2]).
\end{equation*}

3. Solve the following three reduced matrix 
differential equations on the low-dimensional subspaces:
\begin{align*}
	\dot{B}(t)
	&= F\bigl(t,B(t)Q_2^\top\bigr)Q_2,
	&
	B(t_0)&=N_0Q_2,
	\\
	\dot{C}(t)
	&= F\bigl(t,Q_1C(t)^\top\bigr)^\top Q_1,
	&
	C(t_0)&=N_0^\top Q_1,
	\\
	\dot{D}(t)
	&= Q_1^\top F\bigl(t,Q_1D(t)Q_2^\top\bigr)Q_2,
	&
	D(t_0)&=Q_1^\top N_0Q_2.
\end{align*}
These three reduced equations can be solved simultaneously.

4. At $t=t_0+\tau$, assemble the rank-$r$ 
approximation as
\begin{equation*}
	N_1
	=
	B(t_0+\tau)
	\bigl(D(t_0+\tau)\bigr)_r^\dagger
	C(t_0+\tau),
\end{equation*}
where $(\cdot)_r^\dagger$ denotes the truncated pseudo-inverse to rank $r$.

The detailed derivation and analysis of the DGN method can be found 
in \cite{carrel2025randomized}, and the complete procedure is summarized 
in Algorithm \ref{alg:dgn}.

\section{Implementation and extension}
\label{sec3}

In this section, we present the core implementation of our splitting-based RDLR methods, and extend them to the rank-adaptive scenario.

\subsection{Implementation}\label{subsec3.1}

For Eq.~\eqref{eq1.1}, we focus on how to compute the low-rank solution at $t_1 = t_0 + \tau$. The low-rank solutions at other time steps ($t_n = t_0 + n\tau$ for $n \geq 2$) can be obtained analogously. 

For the linear subproblem \eqref{eq2.1}, the core operation involves computing the action of the matrix exponential $\mathrm{e}^{\tau A}$ (Lie-Trotter splitting) or $\mathrm{e}^{\frac{\tau}{2} A}$ (Strang splitting) on the current low-rank solution factor. Although efficient algorithms exist for matrix exponential operations (e.g., Taylor series expansion \cite{al2011computing}, Leja point interpolation \cite{Caliari2016}, and Krylov subspace approximation \cite{guttel2013rational,saad1992analysis}), we utilize MATLAB's built-in \texttt{expm} function due to its robustness and numerical stability. Crucially, the matrix exponential $\mathrm{e}^{\tau A}$ or $\mathrm{e}^{\frac{\tau}{2} A}$ is computed only once, prior to the main time integration loop, and reused throughout the time integration, providing significant computational savings. 

The nonlinear subproblem \eqref{eq2.2} is solved via using the DRSVD method (i.e., Algorithm \ref{alg:drsvd}) or the DGN method (i.e., Algorithm \ref{alg:dgn}). Within this framework, the underlying matrix differential equation is numerically integrated using a fourth-order Runge-Kutta (RK4). To ensure sufficient temporal resolution and accuracy while maintaining computational efficiency, we set $10$ uniform substeps within each main time step of size $\tau$. Moreover, our proposed splitting-based RDLR approximations are denoted as DRSVD-LT, DRSVD-ST, DGN-LT, DGN-ST. For example, DRSVD-LT means the DRSVD method is used in the low-rank Lie-Trotter splitting \eqref{eq:star1}. Other notations are defined in a similar manner. In these notations, the distinguishing labels "LT" and "ST" are added to distinguish which splitting method is choose.

For clarification, the DRSVD-LT and DRSVD-ST methods are summarized in Algorithms \ref{alg:lie_drsvd} and \ref{alg:strang_drsvd}, respectively. Algorithms for the DGN-LT and the DGN-ST methods can be obtained by replacing the DRSVD method with  the DGN method in Algorithms \ref{alg:lie_drsvd} and \ref{alg:strang_drsvd}, respectively.
Next, we will extend these approximations to the rank-adaptive scenario.

\begin{algorithm}
	\caption{DRSVD-LT method for \eqref{eq1.1}, single time step} 
\label{alg:lie_drsvd}
$\texttt{DRSVD-LT}(U, S, V, \tau, A, F, r, p, q)$
\begin{algorithmic}[1]
	\REQUIRE Factors $U, S, V$ of rank-$r$ approximation $X_0 \approx X(t_0)$ with $U \in \mathbb{R}^{m \times r}$, $V \in \mathbb{R}^{n \times r}$, $S \in \mathbb{R}^{r \times r}$, step size $\tau \geq 0 $, matrix $A$, nonlinear function $F$, target rank $r \geq 0 $, oversampling $p \geq 0 $, power iterations $q \geq 0 $
	
	\noindent\textit{\textbf{The nonlinear term \eqref{eq2.2} (full-step) -- apply matrix exponential to linear part}}
	
	\STATE  $U_1, S_1, V_1 = \texttt{DynamicalRandomizedSVD}(U, S, V, F, \tau, r, p, q)$ via Algorithm \ref{alg:drsvd}
	
	\textit{\textbf{The linear term \eqref{eq2.1} (full-step) -- apply matrix exponential to linear part}}
	
	\STATE Compute $QR$-decomposition of ${U}_2 R_U = \mathrm{e}^{\tau A} {U}_1$, $R_U$ invertible, ${U}_2$ orthonormal
	
	\STATE Compute $QR$-decomposition of ${V}_2 R_V = \mathrm{e}^{\tau A} V_1$, $R_V$ invertible, ${V}_2$ orthonormal
	
	\STATE Set ${S}_2 = R_U S_1 R_V^\top$
	
	\RETURN {Factors $U_2, S_2, V_2$ of rank-$r$ approximation $X_1 = U_2 S_2 V_2^\top \approx X(t_0 + \tau)$}
\end{algorithmic}
\end{algorithm}

\begin{algorithm}
\caption{DRSVD-ST method for \eqref{eq1.1}, single time step} 
\label{alg:strang_drsvd}
$\texttt{DRSVD-ST}(U, S, V, \tau, A, F, r, p, q)$
\begin{algorithmic}[1]
\REQUIRE Factors $U, S, V$ of rank-$r$ approximation $X_0 \approx X(t_0)$ with $U \in \mathbb{R}^{m \times r}$, $V \in \mathbb{R}^{n \times r}$, $S \in \mathbb{R}^{r \times r}$, step size $\tau \geq 0 $, matrix $A$, nonlinear function $F$, target rank $r \geq 0 $, oversampling $p \geq 0 $, power iterations $q \geq 0 $

\noindent\textit{\textbf{The linear term \eqref{eq2.1} (half-step) -- apply matrix exponential to linear part}}

\STATE Compute $QR$-decomposition of ${U}_1 R_U = \mathrm{e}^{\frac{\tau}{2} A} {U}$, $R_U$ invertible, ${U}_1$ orthonormal

\STATE Compute $QR$-decomposition of ${V}_1 R_V = \mathrm{e}^{\frac{\tau}{2} A} V$, $R_V$ invertible, ${V}_1$ orthonormal

\STATE Set ${S}_1 = R_U S R_V^\top$

\noindent\textit{\textbf{The nonlinear term \eqref{eq2.2} (full-step) -- solve nonlinear part with DRSVD}}

\STATE $U_2, S_2, V_2 = \texttt{DynamicalRandomizedSVD}(U_1, S_1, V_1, F, \tau, r, p, q)$ via Algorithm \ref{alg:drsvd}

\noindent\textit{\textbf{The linear term \eqref{eq2.1} (half-step) -- apply matrix exponential to linear part}}

\STATE Compute $QR$-decomposition of ${U}_3 R_U^* = \mathrm{e}^{\frac{\tau}{2} A} {U}_2$, $R_U^*$ invertible, ${U}_3$ orthonormal

\STATE Compute $QR$-decomposition of ${V}_3 R_V^* = \mathrm{e}^{\frac{\tau}{2} A} V_2$, $R_V^*$ invertible, ${V}_3$ orthonormal

\STATE Set ${S}_3 = R_U^* S_2 {R_V^*}^\top$

\RETURN {Factors $U_3, S_3, V_3$ of rank-$r$ approximation $X_1 = U_3 S_3 V_3^\top \approx X(t_0 + \tau)$}
\end{algorithmic}
\end{algorithm}

\subsection{Rank-adaptive extension of the proposed methods}\label{subsec3.2}
In most applications, the target rank is typically unknown in advance \cite{zhao2024adaptive}, or varies greatly over time. Consequently, the adaptive randomized rangefinder based on tolerance to dynamically adjust the target rank-$r$ is proposed in \cite{halko2011finding}. This method relies on the following lemma.
\begin{lemma}(\text{\cite{woolfe2008fast}})
Let $A\in\mathbb{R}^{m\times n}$, fix a positive integer $\kappa$ and a real number $\alpha>1$. Draw an independent family $\{\omega^{(i)}:i=1,2,\ldots,\kappa\}$ of standard Gaussian vectors. Then, 
\begin{equation*}
\left\|A\right\|_2\leq\alpha\sqrt{\frac{2}{\pi}}\cdot\max_{i=1,\ldots,\kappa}\left\|A\omega^{(i)}\right\|_2
\end{equation*}
with probability $1-\alpha^{-\kappa}$. In which, $||\cdot||_2$ is the vector 2-norm.
\end{lemma}

Let $A=N(t_0+\tau)-QQ^{\top}N(t_0+\tau)$ and $\alpha=10$ the following conclusion is obtained:
\begin{align*}
&\|N(t_0+\tau)-QQ^{\top}N(t_0+\tau)\|_2 \\
&\quad \leq 10\sqrt{\frac{2}{\pi}} \cdot \max_{i=1,\ldots,\kappa} \biggl\|
\begin{aligned} &N(t_0+\tau)\omega^{(i)}
- QQ^{\top}N(t_0+\tau)\omega^{(i)}\biggr\|_2
\end{aligned}
\end{align*}
with probability $1-10^{-\kappa}$.

Recall that the dynamical rangefinder facilitates the approximation $B(t_0+\tau)\approx N(t_0+\tau)\Omega$. Gaussian matrix $\Omega$, which leads to $N(t_0+\tau)\omega^{(i)}=N(t_0+\tau)\Omega_{[:,i]}\approx B(t_0+\tau)_{[:,i]}$. In other words, $B(t_0+\tau)$ can be used as an accurate error estimate for the range at (almost) no extra cost. In addition, \cite{carrel2025randomized} opts for block computation (adding $\kappa$ columns each time) to cut iterations, 
%
%
%
%
%
%
balancing accuracy and speed. This procedure can be implemented by Algorithm \ref{alg:adaptive_rangefinder}. This approach offers a flexible and efficient adaptive framework, making it highly suitable for large-scale problems with unknown or time-dependent target ranks.

Similar to the duality transformation from the dynamical rangefinder to the dynamical corangefinder, the adaptive dynamical corangefinder is straightforward. Specifically, by replacing the B-step in the adaptive dynamical rangefinder (Algorithm \ref{alg:adaptive_rangefinder}) with its dual C-step formulation \eqref{eq2.7}.

Within the rank-adaptive framework, we systematically extend DRSVD and DGN to the rank-adaptive scenario. For the adaptive DRSVD, the primary modification involves replacing fixed-rank truncation with tolerance-based truncation. This eliminates the need for prespecified oversampling parameters while preserving the rank-preserving properties for the stiff term integration. Analogously, the adaptive DGN implements rank optimization by replacing both rangefinder and corangefinder steps with adaptive processes. Moreover, our methods under the rank-adaptive framework are denoted as ADRSVD-LT, ADRSVD-ST, ADGN-LT and ADGN-ST. 

We prioritize this strategy as it not only removes subjectivity in oversampling parameter selection but also permits greater rank growth flexibility to capture transient solution features, enhancing adaptability for high-dimensional dynamical systems. The numerical results of the proposed methods are presented below.

\section{Numerical results}
\label{sec4}

This section validates the performance of the proposed algorithm through numerical experiments. In Section \ref{sec4.1}, we quantitatively evaluate the convergence order and computational accuracy using benchmark cases: the Allen-Cahn equation and the differential Riccati equation (DRE), with comparative analysis against existing splitting low-rank techniques. Section \ref{sec4.2} examines computational efficiency and stability in complex scenarios via  long-time simulations, providing theoretical foundations for practical implementations. 

\subsection{Convergence tests}\label{sec4.1} 

This section aims to validate the convergence behavior of the proposed DRLR methods on two typical stiff problems. In the following, we will conduct a comparative analysis against the low-rank techniques in \cite{ostermann2019convergence,lubich2014projector,kieri2019projection,ceruti2022rank,lam2025randomized,carrel2023projected,ceruti2022unconventional}. By varying time steps and spatial grid points, we evaluate convergence orders and computational accuracy of our methods, while presenting essential quantitative metrics for each case.

In these examples, for convenience, we set the number of grid points in both the $x$- and \(y\)-directions to $\hat{N}(\hat{N}\in \mathbb{Z}^+)$ and the grid spacing to $h$. We define the relative error at time $T = t_0 + {\hat{M}}\tau$ as:
\begin{equation*}
\text{relerr}(\tau, h) = \frac{\left\| X_{\hat{M}} - X(T) \right\|_{F}}{\left\| X(T) \right\|_{F}},
\end{equation*}
where $X_{\hat{M}}$ is the numerical approximation of $X(T)$. Then, we denote
\begin{equation*}
\text{rate}_{\tau} = \log_{\tau_1/\tau_2} \frac{\text{relerr}(\tau_1,h)}{\text{relerr}(\tau_2,h)} 
\quad \text{and} \quad
\text{rate}_{h} = \log_{h_1/h_2} \frac{\text{relerr}(\tau,h_1)}{\text{relerr}(\tau,h_2)}.
\end{equation*}

\subsubsection{Allen-Cahn equation}\label{sec4.1.1} 
The Allen-Cahn equation, a standard category of stiff matrix differential equations, serves as a classical phase-field modeling framework predominantly employed to elucidate phase separation, interface dynamics, and the evolution of multiphase systems \cite{allen1972ground,allen1973correction}. The most basic representation of this equation includes a nonlinear term in the form of a cubic polynomial (reaction term) alongside a diffusion term \cite{kassam2005fourth}, specifically:
\begin{equation}\label{eq4.1}
\frac{\partial f}{\partial t}=\varepsilon\Delta f+f-f^3,
\end{equation}
in this context, we set $\varepsilon=0.01$. We designate the initial value as follows:
\begin{equation*}
f_0(x,y)=\frac{\left[\mathrm{e}^{-\tan^2(x)}+\mathrm{e}^{-\tan^2(y)}\right]\sin(x)\sin(y)}{1+\mathrm{e}^{|\csc(-x/2)|}+\mathrm{e}^{|\csc(-y/2)|}}.
\end{equation*}
\begin{table}[ht]
\tabcolsep=3.0pt
\renewcommand{\arraystretch}{1.3} 
\begin{center}
\caption{Relative errors of existing low-rank approximation methods at $r=12$ after two temporal steps ($N=1024, T=0.1$), reference solution truncated to target rank.}
\centering
\scriptsize
\begin{tabular}{lr}
	\hline\noalign{\hrule height 1pt}
	Method & Relative error after two steps \\
	\cline{1-1}\cline{2-2}
	Projector-splitting \cite{lubich2014projector} & $\dagger$ \\
	Projected RK1 (PRK1) \cite{kieri2019projection}  & $2.068177\text{e-}03$ \\
	Randomized low-rank RK1 (RRK1) \cite{lam2025randomized}  & $2.068003\text{e-}03$ \\
	Projected exponential RK1 (PERK1) \cite{carrel2023projected}  & $2.216035\text{e-}03$ \\
	Projected RK2 (PRK2) \cite{kieri2019projection} & $3.260517\text{e-}05$ \\
	Randomized low-rank RK2 (RRK2) \cite{lam2025randomized} & $3.196629\text{e-}05$ \\
	Projected exponential RK2 (PERK2) \cite{carrel2023projected} & $3.809570\text{e-}05$ \\
	BUG integrator \cite{ceruti2022unconventional} & $\dagger$ \\
	Augmented BUG integrator \cite{ceruti2022rank}  & $\dagger$ \\
	Best rank-$12$ reference truncation & $8.083527\text{e-}07$ \\
	DRSVD-LT & $1.344408\text{e-}06$ \\
	DRSVD-ST & $1.173795\text{e-}06$ \\
	DGN-LT & $1.344408\text{e-}06$ \\
	DGN-ST & $1.124020\text{e-}06$ \\
	\hline
\end{tabular}
\label{tab_new}
\end{center}
\end{table}

The governing equations are numerically solved on the computational domain $\Omega=[0,2\pi]^{2}$, which is discretized using a uniformly spaced grid with $1024\times1024$ nodal points. Periodic boundary conditions are imposed along all domain boundaries. Through implementation of a second-order centered finite difference discretization, Eq.~\eqref{eq4.1} is transformed into a system of ODEs. Compared to Eq.~\eqref{eq1.1}, the matrix \( A \) becomes the corresponding discrete Laplace matrix, and $A = A^\top$. The nonlinear term \( F(t, X(t)) \) is replaced by \( X(t) - X(t)^{*3} \), where $^{*3}$ is the cubic operation on elementwise. More precisely, 
\begin{equation*}
\dot{X}(t)=AX(t)+X(t)A+X(t)-X(t)^{*3}.
\end{equation*}

Before systematically investigating convergence orders, we first evaluate the robustness of various methods under a large step size. In Table \ref{tab_new}, we compare the performance of our methods (more precisely, DRSVD-LT, DRSVD-ST, DGN-LT and DGN-ST) with some existing DLRA methods applied to \eqref{eq4.1}, where $\dagger$ represents that the method did not converge. The spatial grid size is set to \( N = 1024 \), the final time to \( T = 0.1 \), and only two time steps are taken, corresponding to a step size of \( h = 0.05 \). The target rank is set to \( r = 12 \). The results show that existing DLRA methods perform poorly due to the stiffness of the problem and the large time step. In contrast, the proposed four methods (i.e., DRSVD-LT, DRSVD-ST, DGN-LT and DGN-ST) methods remain computationally stable and achieve higher accuracy.
\begin{figure}[ht]	
\centering
\subfigure{
\includegraphics[width=2.08in,height=2.45in]{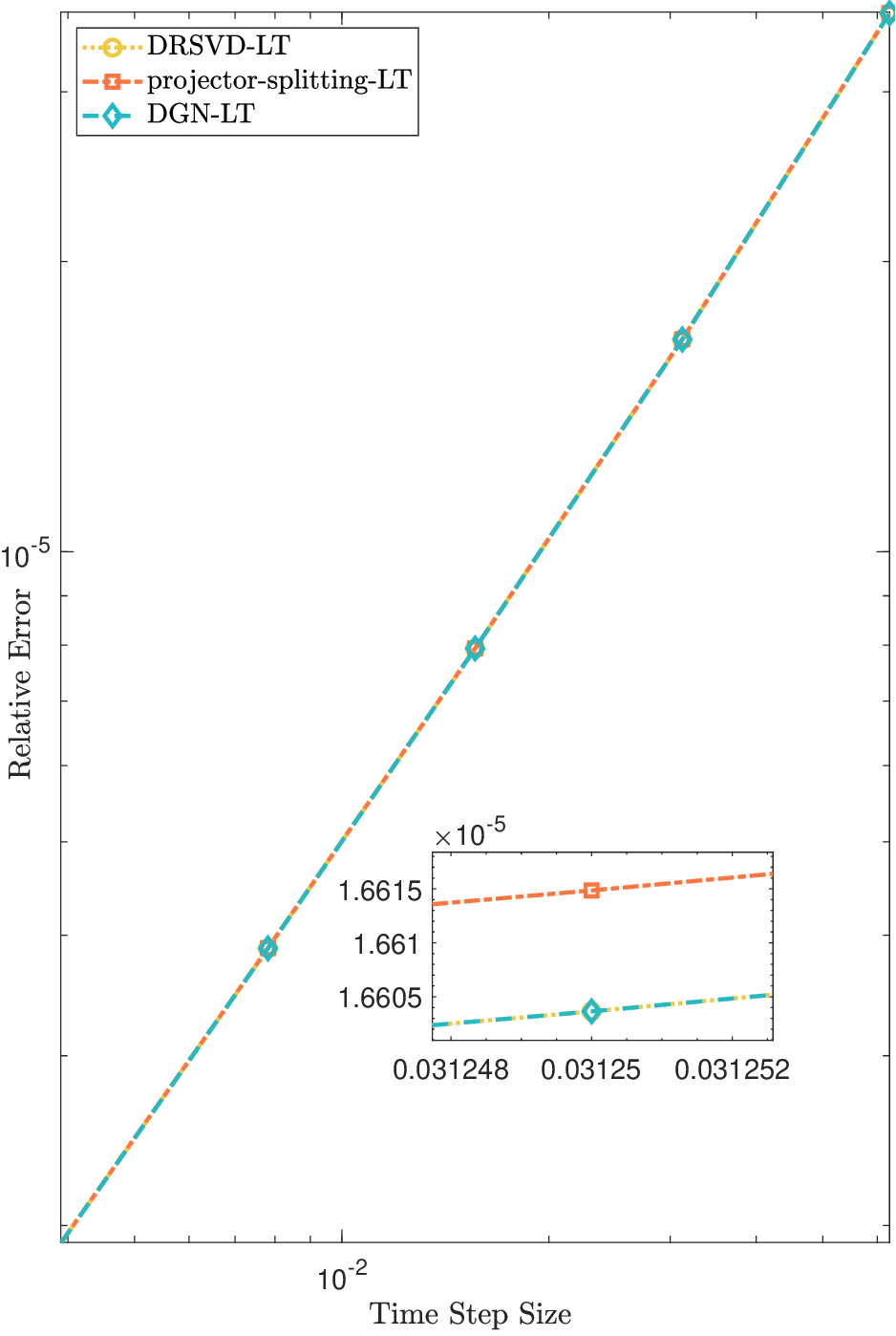}} \hspace{10mm}
\subfigure{
\includegraphics[width=2.08in,height=2.45in]{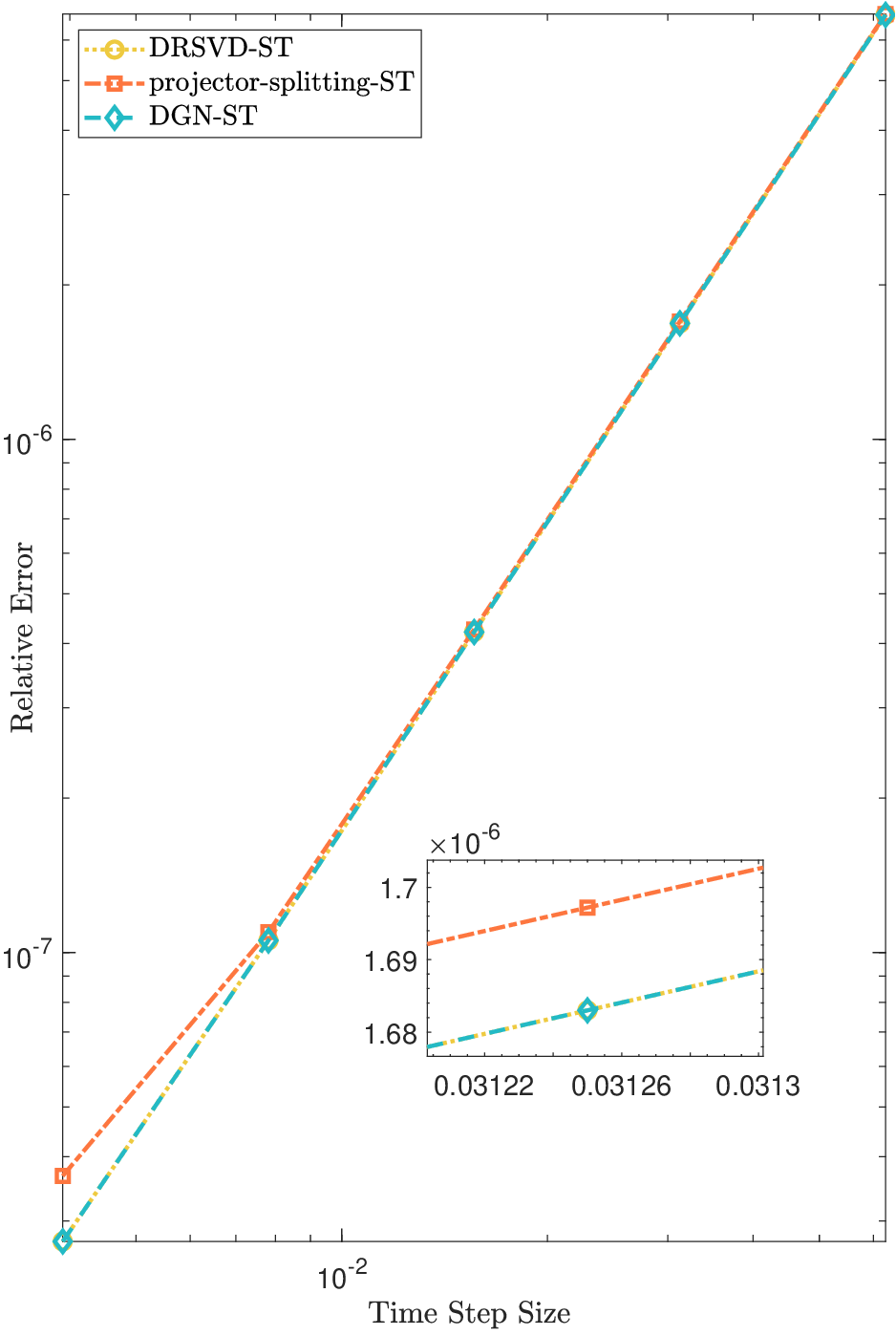}}			
\caption{Comparison of relative errors between DRSVD-type, DGN-type splittings and projector-splitting with different time step size $\tau$ for Eq.~\eqref{eq4.1} for target rank $r = 16$ and the number of time steps $\hat{M} = 1024$.
Left: The Lie-Trotter splitting is employed for decoupling the equation.
Right: The Strang splitting is employed for decoupling the equation.}
\label{fig1}
\end{figure}

\begin{figure}[H]	
\centering
\subfigure{
\includegraphics[width=2.08in,height=2.45in]{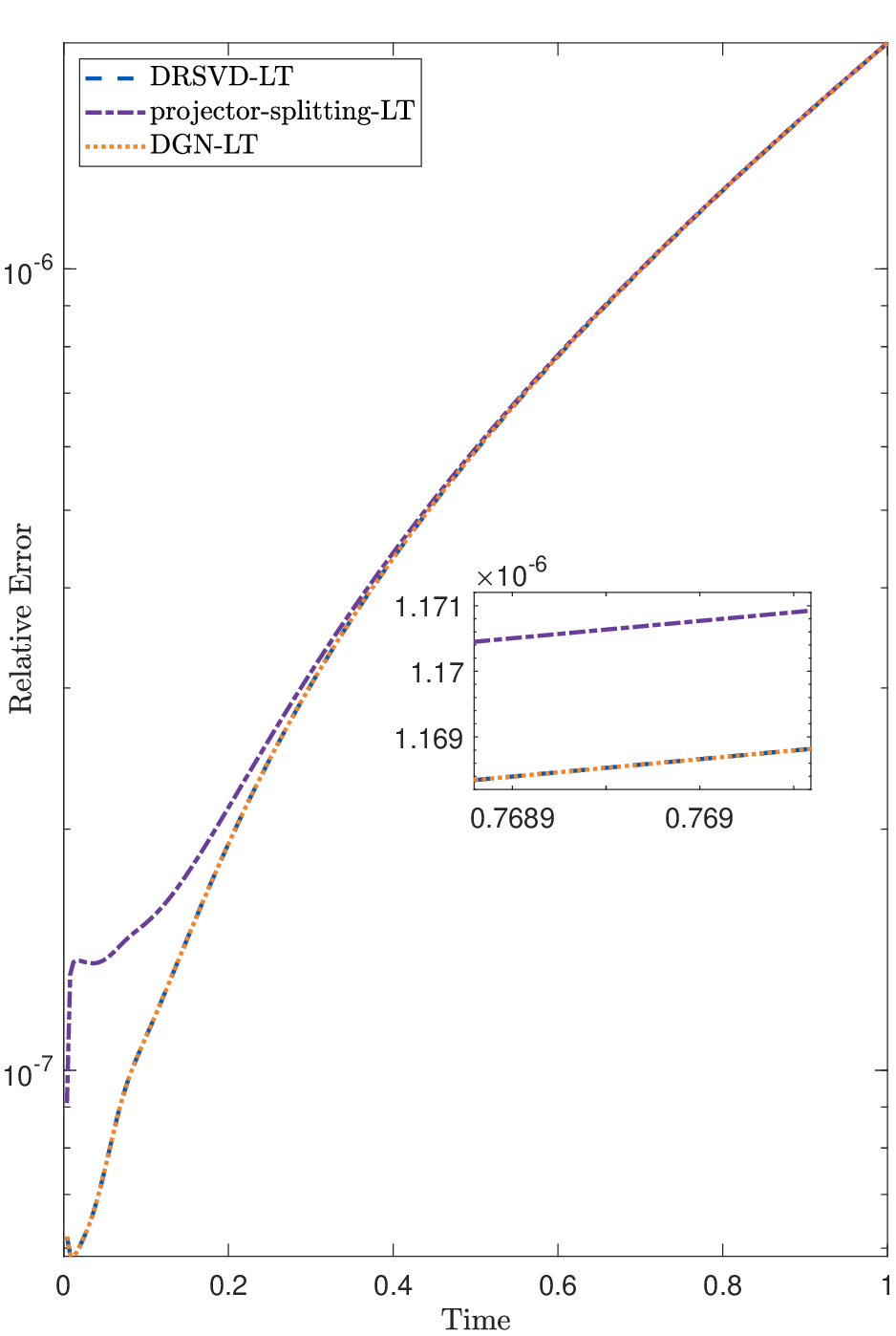}} \hspace{10mm}
\subfigure{
\includegraphics[width=2.08in,height=2.45in]{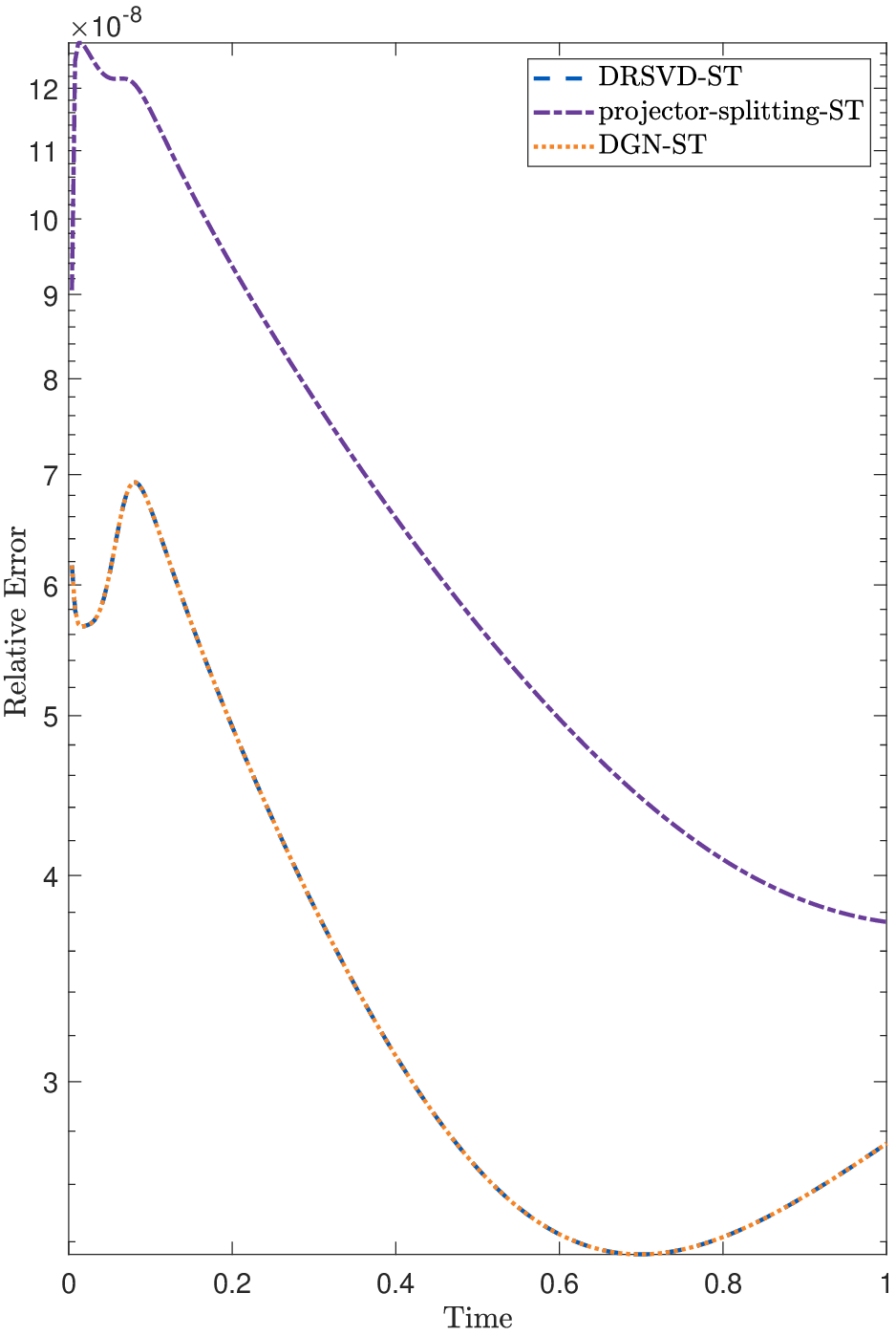}}			
\caption{Comparison of relative errors between DRSVD-type, DGN-type splittings and projector-splitting over time for Eq.~\eqref{eq4.1} for target rank $r = 16$ and the number of time steps $\hat{M} = 1024$.
Left: The  Lie-Trotter splitting is employed for decoupling the  equation.
Right: The Strang splitting is employed for decoupling the  equation.}
\label{fig2}
\end{figure}

Having confirmed the robustness of the proposed methods under a large step size, we now systematically investigate their convergence behavior. The time interval is set to $[0, T]$ with $T=1$. Fig.~\ref{fig1} demonstrates that Eq.~\eqref{eq4.1} solved at low rank $r = 16$, DRSVD-LT, DRSVD-ST, DGN-LT and DGN-ST consistently achieve lower relative errors across various time step sizes $\tau$ compared to the conventional projector-splitting method in \cite{ostermann2019convergence}. This precision advantage holds for both Lie-Trotter and Strang splitting schemes. Fig.~\ref{fig2} further reveals that the four proposed methods exhibit superior stability, characterized by weaker error accumulation over time. Throughout the simulation, the error curves for both proposed methods remain below that of projector-splitting. This stability advantage is especially evident when using the Strang splitting scheme.

Tables \ref{tab1} and \ref{tab2} present the relative errors in Frobenius norm and temporal convergence orders for Eq.~\eqref{eq4.1} $(\hat{N} = 1024)$ using DRSVD-LT, DGN-LT, DRSVD-ST and DGN-ST methods.
For Lie-Trotter splitting (Table \ref{tab1}), both DRSVD-LT and DGN-LT exhibit the first-order convergence rate with negligible differences in errors. On the other hand, relative errors decrease consistently with refining $\hat{M}$. 
DRSVD-ST and DGN-ST (Table \ref{tab2}) achieves near second-order convergence for $r\geq16$. While DRSVD-type splitting or DGN-type splitting yield nearly identical results, minor discrepancies in convergence rates suggest slight method-dependent variations. Overall, Strang splitting (DRSVD-ST and DGN-ST) demonstrates superior accuracy and theoretical second-order behavior, whereas Lie-Trotter splitting (DRSVD-LT and DGN-LT) maintains robust first-order convergence.

\begin{table}[H]\tabcolsep=3.0pt
\renewcommand{\arraystretch}{1.3} 
\begin{center}
\caption{Relative errors in Frobenius norm and temporal convergence orders of DRSVD-LT and DGN-LT ($\hat{N}=1024$) for Eq.~\eqref{eq4.1}.}
\centering
\scriptsize
\begin{tabular}{cccccccccc}
	\hline
	& & \multicolumn{2}{c}{$r = 12$} & \multicolumn{2}{c}{$r = 14$} & \multicolumn{2}{c}{$r = 16$} 
	& \multicolumn{2}{c}{$r = 18$}\\
	\cmidrule(lr){3-4} \cmidrule(lr){5-6} \cmidrule(lr){7-8} \cmidrule(lr){9-10}
	Method & $\hat{M}$ & $\textrm{relerr}(\tau,h)$ & $\textrm{rate}_\tau$ 
	& $\textrm{relerr}(\tau,h)$ & $\textrm{rate}_\tau$
	& $\textrm{relerr}(\tau,h)$ & $\textrm{rate}_\tau$ & $\textrm{relerr}(\tau,h)$ & $\textrm{rate}_\tau$\\
	\hline
	& 16 & 3.6201E-5 & -- & 3.6192E-5 & -- & 3.6192E-5 & -- & 3.6192E-5 & -- \\
	& 32 & 1.6734E-5 & 1.1133 & 1.6604E-5 & 1.1242 & 1.6604E-5 & 1.1242 & 1.6604E-5 & 1.1242 \\
	DRSVD-LT & 64 & 8.5011E-6 & 0.9770 & 7.9359E-6 & 1.0651 & 7.9356E-6 & 1.0651 & 7.9356E-6 & 1.0651 \\
	& 128 & 4.9691E-6 & 0.7746 & 3.9298E-6 & 1.0139 & 3.8773E-6 & 1.0333 & 3.8773E-6 & 1.0333 \\
	& 256 & 3.6753E-6 & 0.4351 & 2.0434E-6 & 0.9435 & 1.9162E-6 & 1.0168 & 1.9162E-6 & 1.0168 \\ 
	\hline 
	& 16 & 3.6201E-5 & -- & 3.6192E-5 & -- & 3.6192E-5 & -- & 3.6192E-5 & -- \\
	& 32 & 1.6734E-5 & 1.1133 & 1.6604E-5 & 1.1242 & 1.6604E-5 & 1.1242 & 1.6604E-5 & 1.1242 \\
	DGN-LT & 64 & 8.4642E-6 & 0.9833 & 7.9359E-6 & 1.0651 & 7.9356E-6 & 1.0510 & 7.9356E-6 & 1.0651 \\
	& 128 & 5.0067E-6 & 0.7575 & 3.9298E-6 & 1.0139 & 3.8773E-6 & 1.0333 & 3.8773E-6 & 1.0333 \\
	& 256 & 3.6997E-6 & 0.4365 & 2.0443E-6 & 0.9429 & 1.9162E-6 & 1.0168 & 1.9162E-6 & 1.0168 \\
	\hline
\end{tabular}
\label{tab1}
\end{center}
\end{table}

\begin{table}[H]\tabcolsep=3.0pt
\renewcommand{\arraystretch}{1.3} 
\begin{center}
\caption{Relative errors in Frobenius norm and temporal convergence orders of DRSVD-ST and DGN-ST ($\hat{N}=1024$) for Eq.~\eqref{eq4.1}.}
\centering
\scriptsize
\begin{tabular}{cccccccccc}
	\hline
	& & \multicolumn{2}{c}{$r = 12$} & \multicolumn{2}{c}{$r = 14$} & \multicolumn{2}{c}{$r = 16$} & \multicolumn{2}{c}{$r = 18$}\\
	\cmidrule(lr){3-4} \cmidrule(lr){5-6} \cmidrule(lr){7-8} \cmidrule(lr){9-10}
	Method &$\hat{M}$ & $\textrm{relerr}(\tau,h)$ & $\textrm{rate}_\tau$ 
	& $\textrm{relerr}(\tau,h)$ & $\textrm{rate}_\tau$
	& $\textrm{relerr}(\tau,h)$ & $\textrm{rate}_\tau$ & $\textrm{relerr}(\tau,h)$ & $\textrm{rate}_\tau$\\
	\hline
	& 16 & 6.7564E-6 & --         & 6.7175E-6 & --         & 6.7175E-6 & --         & 6.7175E-6 & --         \\
	& 32 & 2.6625E-6 & 1.3435     & 1.6834E-6 & 1.9965     & 1.6830E-6 & 1.9969     & 1.6830E-6 & 1.9969     \\
	DRSVD-ST & 64 & 2.8422E-6 & -0.0942    & 4.2343E-7 & 1.9912     & 4.2124E-7 & 1.9983     & 4.2120E-7 & 1.9984     \\
	& 128 & 3.0821E-6 & -0.1169    & 5.8686E-7 & -0.4709    & 1.0558E-7 & 1.9963     & 1.0536E-7 & 1.9991     \\
	& 256 & 3.1086E-6 & -0.0124    & 6.8910E-7 & -0.2317    & 2.7357E-8 & 1.9483     & 2.6387E-8 & 1.9975     \\
	\hline
	& 16 & 6.7564E-6 & --         & 6.7175E-6 & --         & 6.7175E-6 & --         & 6.7175E-6 & --         \\
	& 32 & 2.6625E-6 & 1.3435     & 1.6834E-6 & 1.9965     & 1.6830E-6 & 1.9969     & 1.6830E-6 & 1.9969     \\
	DGN-ST & 64 & 2.9216E-6 & -0.1340     & 4.2343E-7 & 1.9912     & 4.2124E-7 & 1.9983     & 4.2120E-7 & 1.9984     \\
	& 128 & 3.1596E-6 & -0.1130     & 5.9090E-7 & -0.4808    & 1.0558E-7 & 1.9963     & 1.0536E-7 & 1.9991     \\
	& 256 & 3.1594E-6 & 0.0001     & 6.9048E-7 & -0.2247    & 2.7357E-8 & 1.9483     & 2.6387E-8 & 1.9975     \\
	\hline
\end{tabular}
\label{tab2}
\end{center}
\end{table}

\begin{figure}[ht]	
\centering
\subfigure{
\includegraphics[width=2.08in,height=2.45in]{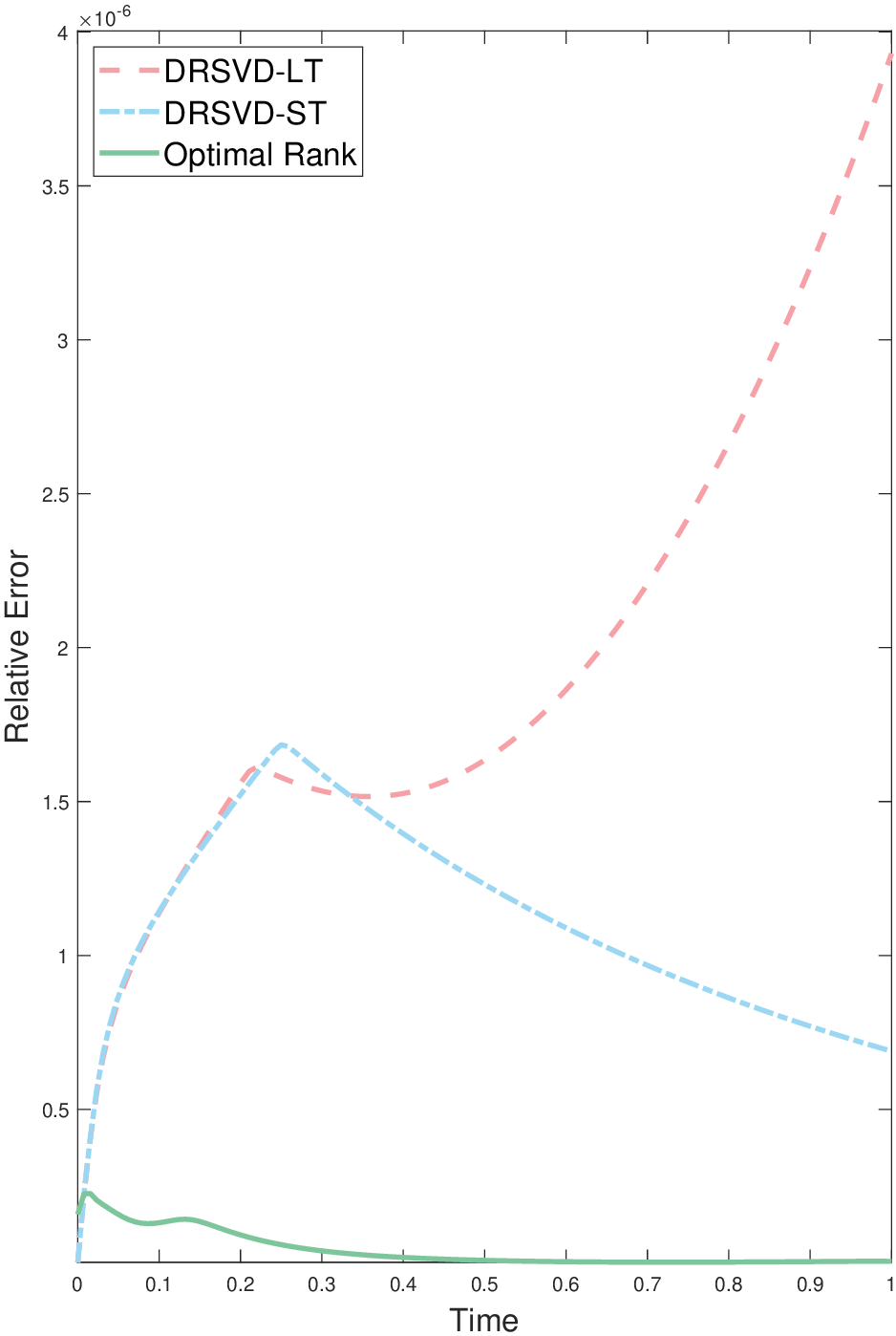}} \hspace{10mm}
\subfigure{
\includegraphics[width=2.08in,height=2.45in]{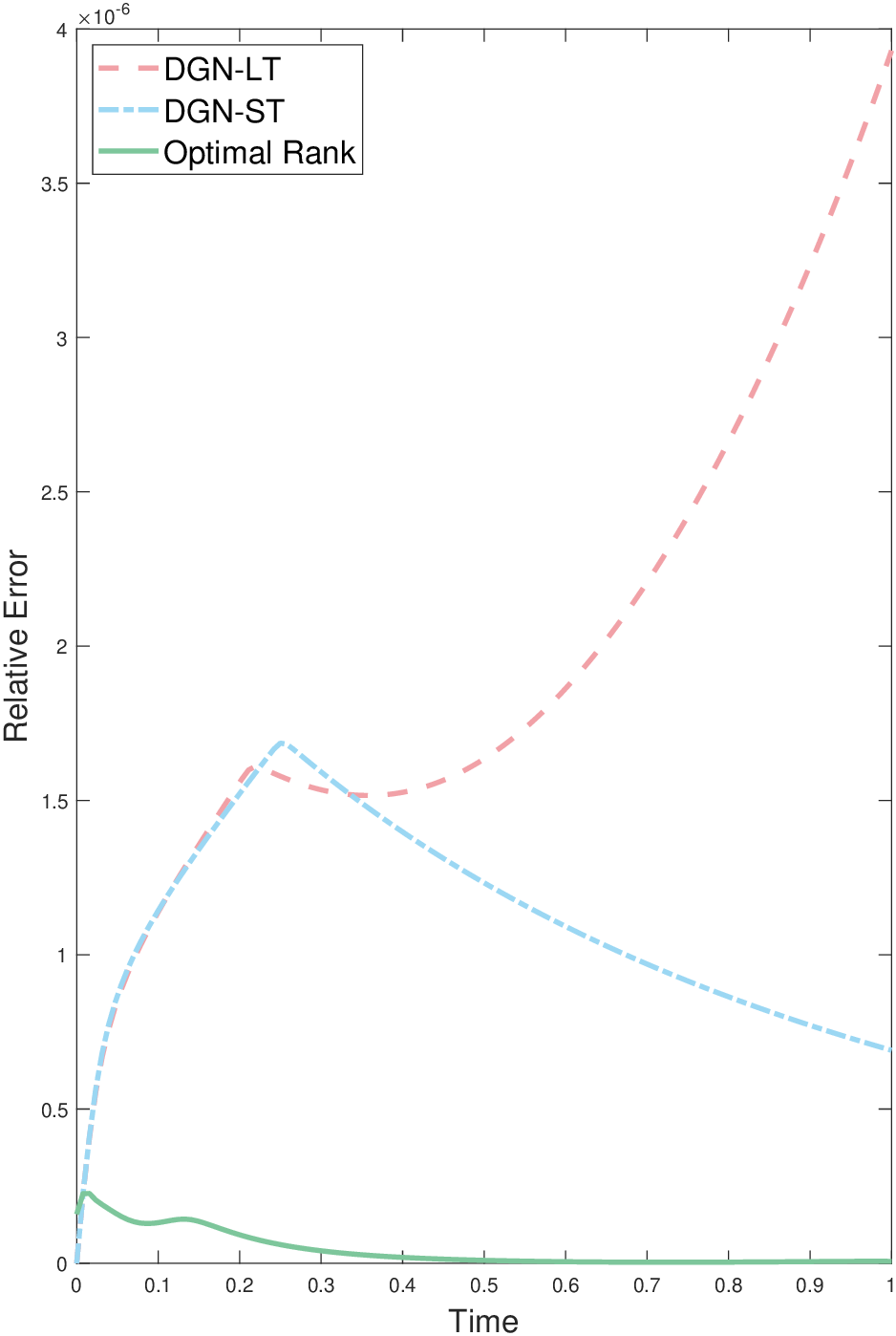}}			
\caption{Comparison of relative errors between Lie-Trotter splitting, Strang splitting and optimal rank approximation over time for Eq.~\eqref{eq4.1} for target rank $r = 14 $ and the number of time steps $\hat{M} = 128$.
Left: The DRSVD method is employed for the nonlinear terms.
Right: The DGN method is employed for the nonlinear terms.}
\label{fig3}
\end{figure}

Fig.~\ref{fig3} illustrates the temporal evolution of errors for both DRSVD-type splitting or DGN-type splitting methods, compared to the optimal rank approximation error, under the conditions of a fixed rank $r = 14$ and the number of time steps $\hat{M} = 128$. The left and right subplots represent the results obtained using the DRSVD and DGN techniques for addressing the nonlinear term, respectively. It is evident from the figure that, regardless of the nonlinear treatment employed (DRSVD or DGN), the error curve of the Strang splitting method consistently remains below that of Lie-Trotter splitting method, thereby confirming its inherent advantage in achieving the second-order accuracy. 
\begin{table}[ht]\tabcolsep=5.0pt
\renewcommand{\arraystretch}{1.3} 
\begin{center}
\caption{Relative errors in Frobenius norm and temporal convergence orders of ADRSVD and ADGN ($\hat{N}=1024$) for Eq.~\eqref{eq4.1}.}
\centering
\scriptsize
\begin{tabular}{ccccccccc}
	\hline
	& \multicolumn{2}{c}{ADRSVD-LT} & \multicolumn{2}{c}{ADGN-LT} & \multicolumn{2}{c}{ADRSVD-ST} & \multicolumn{2}{c}{ADGN-ST}\\
	\cmidrule(lr){2-3}\cmidrule(lr){4-5}\cmidrule(lr){6-7}\cmidrule(lr){8-9}
	$\hat{M}$ & $\textrm{relerr}(\tau,h)$ & $\textrm{rate}_\tau$ 
	& $\textrm{relerr}(\tau,h)$ & $\textrm{rate}_\tau$& $\textrm{relerr}(\tau,h)$ & $\textrm{rate}_\tau$ 
	& $\textrm{relerr}(\tau,h)$ & $\textrm{rate}_\tau$\\
	\hline
	16 & 2.4895E-5 & -- & 2.4895E-5 & -- & 6.7175E-6 & -- & 6.7174E-6 & -- \\
	32 & 1.3757E-5 & 0.8557 & 1.3757E-5 & 0.8557 & 1.6830E-6 & 1.9969 & 1.6830E-6 & 1.9969 \\
	64 & 7.2221E-6 & 0.9296 & 7.2221E-6 & 0.9296 & 4.2134E-7 & 1.998  & 4.2138E-7 & 1.9978 \\
	128 & 3.6988E-6 & 0.9654 & 3.6988E-6 & 0.9654 & 1.0593E-7 & 1.9919 & 1.0593E-07 & 1.9919 \\
	256 & 1.8716E-6 & 0.9828 & 1.8716E-6 & 0.9828 & 2.8555E-8 & 1.8913 & 2.8555E-8 & 1.8913 \\
	\hline
\end{tabular}
\label{tab3}
\end{center}
\end{table}

\begin{figure}[ht]
\centering
\subfigure
{\includegraphics[width=2.1in,height=2.1in]{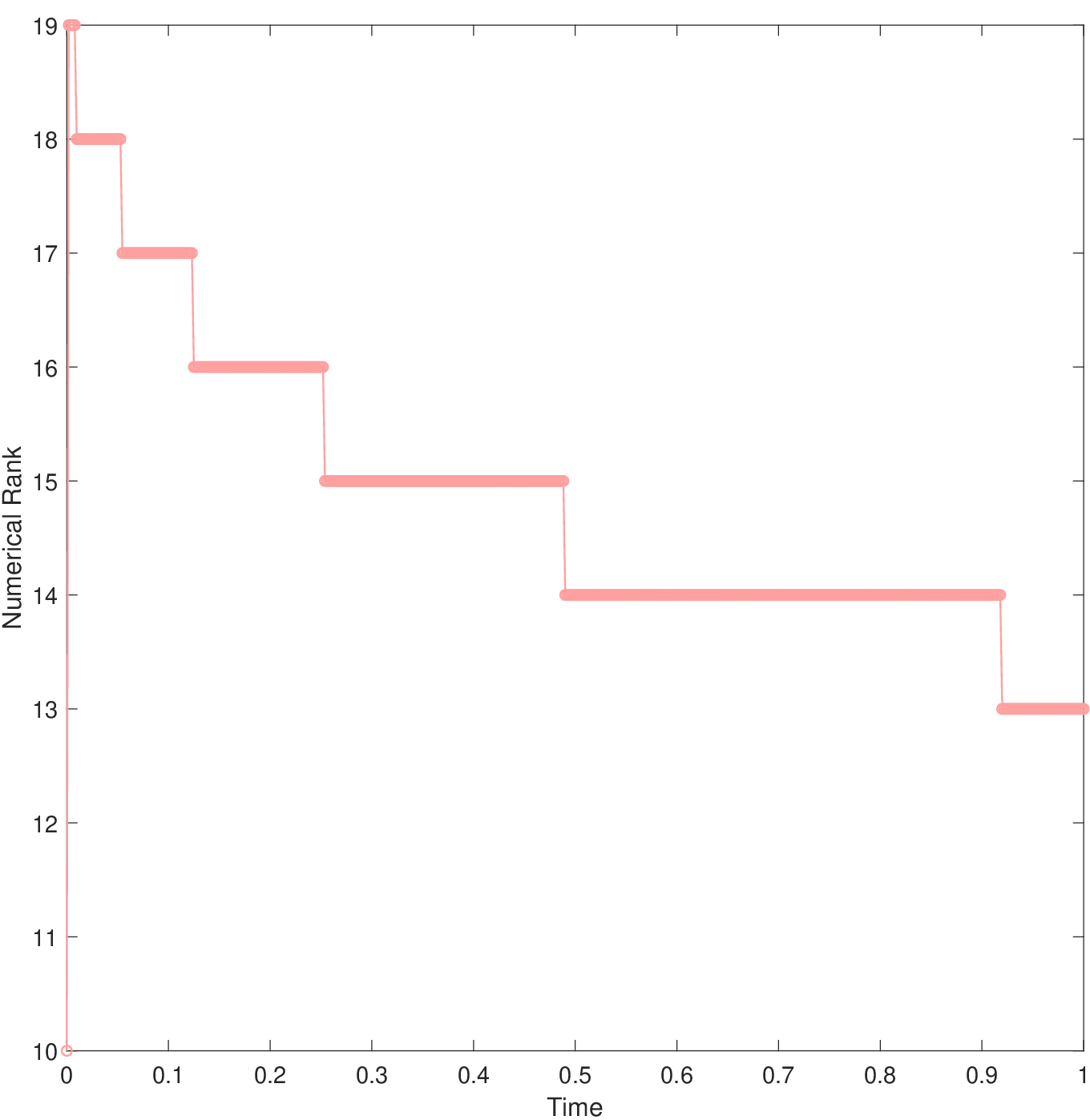}} \hspace{2mm}
\subfigure
{\includegraphics[width=2.1in,height=2.1in]{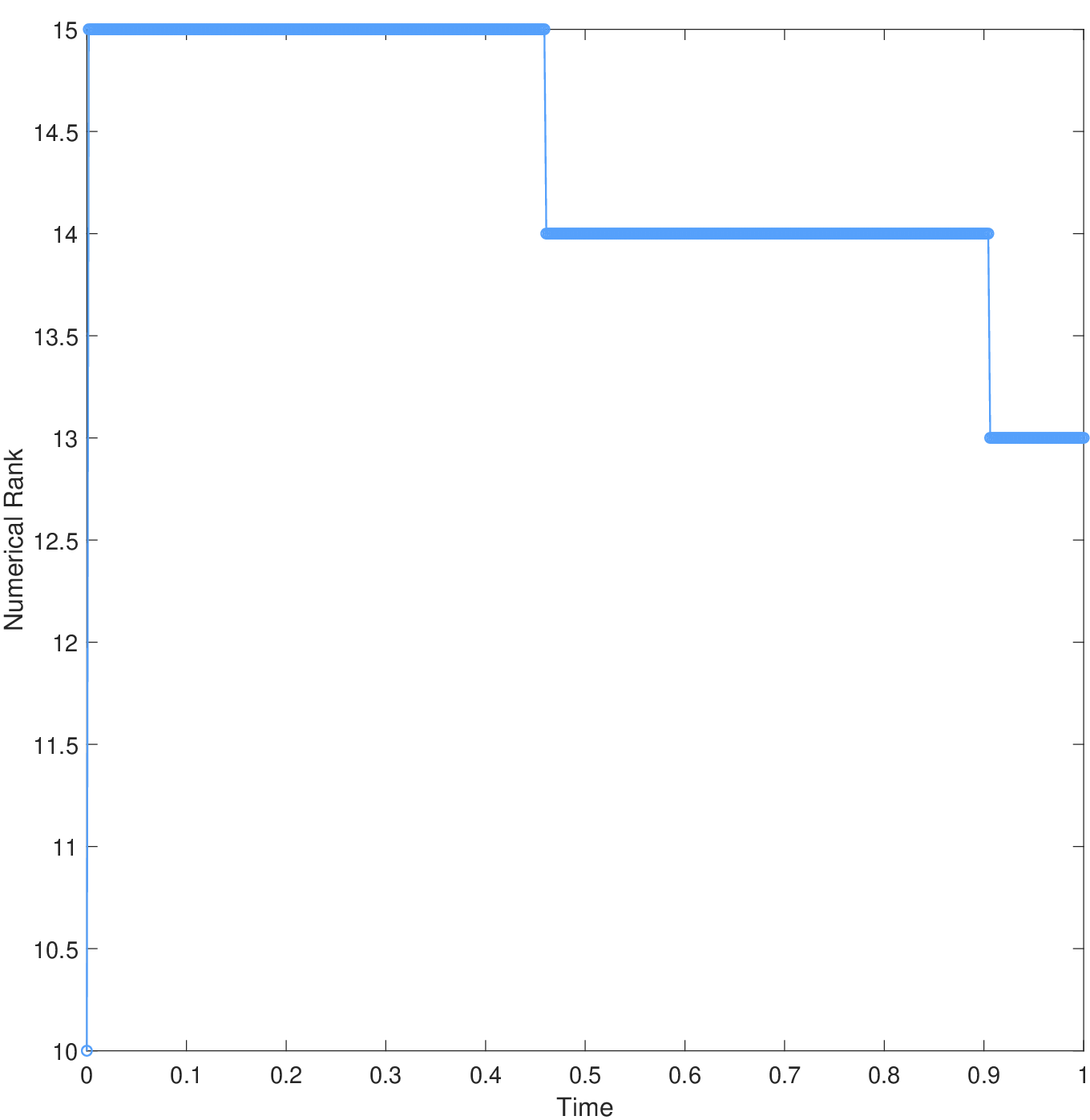}}
\caption{ Numerical rank evolution of Eq.~\eqref{eq4.1} for splitting methods combined with the adaptive low-rank solvers, where $(\hat{M},\hat{N})=(512,1024)$.
Left: ADRSVD-LT; 
Right: ADRSVD-ST.}
\label{fig4}
\end{figure}

Table \ref{tab3} presents temporal convergence behaviors in ADRSVD-LT, ADGN-LT, ADRSVD-ST and ADGN-ST for Eq.~\eqref{eq4.1} with $\hat{N}= 1024$. For the truncation step, the relative tolerance is set to \(\varphi = 10^{-8}\) and the absolute tolerance to \(10^{-12}\). 
The rangefinder step employs a tolerance of \(10^{-8}\). For the Lie-Trotter splitting method, both ADRSVD and ADGN yield identical relative errors in the Frobenius norm, achieving the first-order convergence. In contrast, the Strang splitting method exhibits significantly higher accuracy and approaches the second-order convergence. These results highlight the superior performance of the Strang splitting method and the minimal impact of the nonlinear term approximation method (ADRSVD vs ADGN) on convergence rates. Moreover, Fig.~\ref{fig4} illustrates the evolution of the numerical rank over time when the number of time steps $\hat{M}$ is 256.

Tables \ref{tab4} and \ref{tab5} present the relative errors in the Frobenius norm and the spatial convergence orders for Eq.~\eqref{eq4.1} ($\hat{M} = 1024$) using DRSVD-LT, DGN-LT, DRSVD-ST and DGN-ST methods. For the Lie-Trotter splitting method (Table \ref{tab4}), both DRSVD and DGN demonstrate a spatial convergence behavior where the relative errors decrease as $\hat{N}$ increases. The observed spatial convergence order is approximately 2 for larger number of $\hat{N}$, indicating a consistent and stable convergence pattern. The errors for DRSVD and DGN are quite similar, suggesting comparable performance in capturing the spatial dynamics, with convergence primarily influenced by the spatial discretization process, as expected. For the Strang splitting method (Table \ref{tab5}), a spatial convergence order of around 2 is also achieved, with errors decreasing steadily. The convergence is robust, and the results for DRSVD and DGN are nearly identical, indicating that both methods effectively manage the spatial discretization for Strang splitting.

\begin{table}[H]\tabcolsep=3.0pt
\renewcommand{\arraystretch}{1.3} 
\begin{center}
\caption{Relative errors in Frobenius norm and spatial convergence orders of DRSVD-LT and DGN-LT ($\hat{M} = 1024$) for Eq.~\eqref{eq4.1}.}
\centering
\scriptsize
\begin{tabular}{cccccccccc}
	\hline
	& & \multicolumn{2}{c}{$r = 12$} & \multicolumn{2}{c}{$r = 14$} & \multicolumn{2}{c}{$r = 16$} 
	& \multicolumn{2}{c}{$r = 18$}\\
	\cmidrule(lr){3-4} \cmidrule(lr){5-6} \cmidrule(lr){7-8} \cmidrule(lr){9-10}
	Method & $\hat{N}$ & $\textrm{relerr}(\tau,h)$ & $\textrm{rate}_{h}$ 
	& $\textrm{relerr}(\tau,h)$ & $\textrm{rate}_{h}$
	& $\textrm{relerr}(\tau,h)$ & $\textrm{rate}_{h}$ & $\textrm{relerr}(\tau,h)$ & $\textrm{rate}_{h}$\\
	\hline
	& 16 & 3.5303E-2 & -- & 3.5303E-2 & -- & 3.5303E-2 & -- & 3.5303E-2 & --\\
	& 32 & 8.6882E-3 & 2.0227 & 8.6882E-3 & 2.0227 & 8.6882E-3 & 2.0227 & 8.6882E-3 & 2.0227\\
	DRSVD-LT & 64 & 2.2834E-3 & 1.9279 & 2.2834E-3 & 1.9279 & 2.2834E-3 & 1.9279 & 2.2834E-3 & 1.9279\\
	& 128 & 5.6181E-4 & 2.0230 & 5.6181E-4 & 2.0230 & 5.6181E-4 & 2.0230 & 5.6181E-4 & 2.0230\\
	& 256 & 1.3363E-4 & 2.0718 & 1.3360E-4 & 2.0722 & 1.3359E-4 & 2.0722 & 1.3359E-4 & 2.0722\\
	\hline 
	& 16 & 3.5303E-2 & -- & 3.5303E-2 & -- & 3.5303E-2 & -- & 3.5303E-2 & --\\
	& 32 & 8.6882E-3 & 2.0227 & 8.6882E-3 & 2.0227 & 8.6882E-3 & 2.0227 & 8.6882E-3 & 2.0227\\
	DGN-LT & 64 & 2.2834E-3 & 1.9279 & 2.2834E-3 & 1.9279 & 2.2834E-3 & 1.9279 & 2.2834E-3 & 1.9279\\
	& 128 & 5.6186E-4 & 2.0229 & 5.6185E-4 & 2.0229 & 5.6185E-4 & 2.0229 & 5.6185E-4 & 2.0229\\
	& 256 & 1.3368E-4 & 2.0715 & 1.3364E-4 & 2.0718 & 1.3364E-4 & 2.0718 & 1.3364E-4 & 2.0718\\
	\hline
\end{tabular}
\label{tab4}
\end{center}
\end{table}

\begin{table}[ht]\tabcolsep=3pt
\renewcommand{\arraystretch}{1.3} 
\begin{center}
\caption{Relative errors in Frobenius norm and spatial convergence orders of DRSVD-ST and DGN-ST ($\hat{M} = 1024$) for Eq.~\eqref{eq4.1}.}
\centering
\scriptsize
\begin{tabular}{cccccccccc}
	\hline
	& & \multicolumn{2}{c}{$r = 12$} & \multicolumn{2}{c}{$r = 14$} & \multicolumn{2}{c}{$r = 16$} 
	& \multicolumn{2}{c}{$r = 18$}\\
	\cmidrule(lr){3-4} \cmidrule(lr){5-6} \cmidrule(lr){7-8} \cmidrule(lr){9-10}
	Method & $\hat{N}$ & $\textrm{relerr}(\tau,h)$ & $\textrm{rate}_{h}$ 
	& $\textrm{relerr}(\tau,h)$ & $\textrm{rate}_{h}$
	& $\textrm{relerr}(\tau,h)$ & $\textrm{rate}_{h}$ & $\textrm{relerr}(\tau,h)$ & $\textrm{rate}_{h}$\\
	\hline
	& 16 & 3.5303E-2 & -- & 3.5303E-2 & -- & 3.5303E-2 & -- & 3.5303E-2 & --\\
	& 32 & 8.6882E-3 & 2.0227 & 8.6882E-3 & 2.0227 & 8.6882E-3 & 2.0227 & 8.6882E-3 & 2.0227\\
	DRSVD-ST & 64 & 2.2834E-3 & 1.9279 & 2.2834E-3 & 1.9279 & 2.2834E-3 & 1.9279 & 2.2834E-3 & 1.9279\\
	& 128 & 5.6181E-4 & 2.0229 & 5.6185E-4 & 2.0229 & 5.6185E-4 & 2.0229 & 5.6185E-4 & 2.0229\\
	& 256 & 1.3368E-4 & 2.0715 & 1.3640E-4 & 2.0718 & 1.3364E-4 & 2.0718 & 1.3364E-4 & 2.0718\\
	\hline 
	& 16 & 3.5303E-2 & -- & 3.5303E-2 & -- & 3.5303E-2 & -- & 3.5303E-2 & --\\
	& 32 & 8.6882E-3 & 2.0227 & 8.6882E-3 & 2.0227 & 8.6882E-3 & 2.0227 & 8.6882E-3 & 2.0227\\
	DGN-ST & 64 & 2.2834E-3 & 1.9279 & 2.2834E-3 & 1.9279 & 2.2834E-3 & 1.9279 & 2.2834E-3 & 1.9279\\
	& 128 & 5.6186E-4 & 2.0229 & 5.6185E-4 & 2.0229 & 5.6185E-4 & 2.0229 & 5.6185E-4 & 2.0229\\
	& 256 & 1.3368E-4 & 2.0715 & 1.3364E-4 & 2.0718 & 1.3364E-4 & 2.0718 & 1.3364E-4 & 2.0718\\
	\hline
\end{tabular}
\label{tab5}
\end{center}
\end{table}

\begin{table}[ht]\tabcolsep=3.0pt
\renewcommand{\arraystretch}{1.3} 
\begin{center}
\caption{Relative errors in Frobenius norm and spatial convergence orders of ADRSVD and ADGN ($\hat{M} = 1024$) for Eq.~\eqref{eq4.1}.}
\centering
\scriptsize
\begin{tabular}{ccccccccc}
	\hline
	& \multicolumn{2}{c}{ADRSVD-LT} & \multicolumn{2}{c}{ADGN-LT} & \multicolumn{2}{c}{ADRSVD-ST} & \multicolumn{2}{c}{ADGN-ST}\\
	\cmidrule(lr){2-3}\cmidrule(lr){4-5}\cmidrule(lr){6-7}\cmidrule(lr){8-9}
	$\hat{N}$ & $\textrm{relerr}(\tau,h)$ & $\textrm{rate}_{h}$ 
	& $\textrm{relerr}(\tau,h)$ & $\textrm{rate}_{h}$& $\textrm{relerr}(\tau,h)$ & $\textrm{rate}_{h}$ 
	& $\textrm{relerr}(\tau,h)$ & $\textrm{rate}_{h}$\\
	\hline
	16 & 3.5303E-2 & -- & 3.5303E-2 & -- & 3.5303E-2 & -- & 3.5303E-2 & -- \\
	32 & 8.6882E-3 & 2.0227 & 8.6882E-3 & 2.0227 & 8.6882E-3 & 2.0227 & 8.6882E-3 & 2.0227 \\
	64 & 2.2834E-3 & 1.9279 & 2.2834E-3 & 1.9279 & 2.2834E-3 & 1.9279 & 2.2834E-3 & 1.9279 \\
	128 & 5.6181E-4 & 2.0230 & 5.6181E-4 & 2.0230 & 5.6185E-4 & 2.0229 & 5.6185E-4 & 2.0229 \\
	256 & 1.3359E-4 & 2.0722 & 1.3359E-4 & 2.0722 & 1.3364E-4 & 2.0718 & 1.3364E-4 & 2.0718 \\
	\hline
\end{tabular}
\label{tab6}
\end{center}
\end{table}

Table~\ref{tab6} illustrates that for Eq.~\eqref{eq4.1}, both the Lie-Trotter splitting and Strang splitting methods, when coupled with ADRSVD and ADGN discretizations, exhibit the second-order spatial convergence. As the spatial grid numbers \(\hat{N}\) increases from $16$ to $256$, the relative errors systematically decrease with increasing \(\hat{N}\), and the convergence rates stabilize near the theoretical value of 2. The errors and convergence rates of ADRSVD and ADGN are indistinguishable within numerical precision, demonstrating equivalent performance in spatial discretization. Preliminary analysis suggests that minor fluctuations in convergence rates at coarse grids (\(\hat{N} \leq 64\)) are primarily due to the dominance of truncation errors, while the predominance of spatial discretization errors at finer grids ensures the recovery of theoretical convergence rates. These results validate the robustness of the methods in addressing nonlinear terms and confirm the theoretical accuracy of the spatial discretization schemes.

\subsubsection{Differential Riccati equation}\label{sec4.1.2} 
Another example of a stiff PDE is called DRE, which plays a crucial role in control theory, filtering and estimation, and financial mathematics by solving problems such as optimal control, state estimation, and portfolio optimization \cite{abou2012matrix,anderson2007optimal,zhou2000continuous}.

In this part, we focus on DRE in the optimal control problem of linear quadratic regulators with parabolic PDEs over a finite time interval $T$. Thus, we construct the following linear control system:
\begin{equation}\label{eq4.2}
\dot{x}=Ax+Bu,\quad x(t_0)=x_0,
\end{equation}
where $A\in\mathbb{R}^{d\times d}$ and $B\in\mathbb{R}^{d\times m}$ are the system matrices, $x\in\mathbb{R}^m$ is the state variable, $u\in\mathbb{R}^m$ is the control. The function $\mathcal{J}$ that need to be minimized is as follows: 
\begin{equation*}
\mathcal{J}(u,x)=\frac{1}{2}\int_{t_0}^{T}\left(x(t)^\top C^\top QCx(t)+u(t)^\top Ru(t)\right)\mathrm{d}t,
\end{equation*}
where $C\in\mathbb{R}^{q\times d}$, ${Q}\in\mathbb{R}^{q\times q}$ is symmetric and positive semidefinite and ${R}\in\mathbb{R}^{m\times m}$ is symmetric and positive definite. Moreover, the optimal control is given in feedback form by $u_{\mathrm{opt}}(t)=-R^{-1}B^{\top}X(t)x(t)$, where $X(t)$ is the solution of the following DRE:
\begin{equation*}
\dot{X}(t)=A^\top X(t)+X(t)A+C^\top QC-X(t)BR^{-1}B^\top X(t).
\end{equation*}

To further demonstrate the applicability of our approaches, we utilize a test case based on the diffusion-convection equation, which was developed from \cite{penzl2000matlab}. Consider the following equation:
\begin{equation}\label{eq4.3}
\partial_tw=\Delta w-10x\partial_xw-100y\partial_yw,\quad w|_{\partial\Omega}=0,
\end{equation}
defined on the domain $\Omega = (0, 1)^2$, subject to homogeneous Dirichlet boundary conditions. The matrix \( A \) arises from standard centered finite difference spatial discretization of (\(\ref{eq4.3}\)) with $\widetilde{d}$ uniformly spaced stiff points per dimension.  Let \(x_{i}=i\delta\) (\(\delta = (\widetilde{d} + 1)^{-1}\), \(i = 1,\dots,\widetilde{d}\)) denote interior discretization points in the $x$-direction. For matrix \(B \in \mathbb{R}^{\widetilde{d} \times 1}\), its entries follow:
\begin{equation*}B_i=
\begin{cases}
1 & \mathrm{if}\quad0.1<x_i\leq0.3, \\
0 & \mathrm{otherwise}.
\end{cases}\end{equation*}
For 
\(C \in \mathbb{R}^{1 \times \widetilde{d}}\), the entries are: 
\begin{equation*}C_i=
\begin{cases}
1 & \mathrm{if}\quad0.7<x_i\leq0.9, \\
0 & \mathrm{otherwise}.
\end{cases}\end{equation*}
Additionally, \(R = I\) and \(Q = 100I\) are selected. Similarly for the $y$-direction.

In the following numerical experiments, we take the initial value \(X_0 = I\), the final time \(T = 0.1\), and the reference solution calculated with the Dormand-Prince 5(4) method (denoted as dop54) \cite{dormand1980family}. We show the results for \(\hat{N} = 400\) grid points of the spatial discretization.

\begin{table}[H]\tabcolsep=3.0pt
\renewcommand{\arraystretch}{1.3} 
\begin{center}
\caption{Relative errors in Frobenius norm and temporal convergence orders of DRSVD-LT and DGN-LT ($\hat{N}=400$) for Eq.~\eqref{eq4.2}.}
\centering
\scriptsize
\begin{tabular}{cccccccccc}
	\hline
	& & \multicolumn{2}{c}{$r = 4$} & \multicolumn{2}{c}{$r = 6$} & \multicolumn{2}{c}{$r = 8$} 
	& \multicolumn{2}{c}{$r = 10$}\\
	\cmidrule(lr){3-4} \cmidrule(lr){5-6} \cmidrule(lr){7-8} \cmidrule(lr){9-10}
	Method & $\hat{M}$ & $\textrm{relerr}(\tau,h)$ & $\textrm{rate}_\tau$ 
	& $\textrm{relerr}(\tau,h)$ & $\textrm{rate}_\tau$
	& $\textrm{relerr}(\tau,h)$ & $\textrm{rate}_\tau$ & $\textrm{relerr}(\tau,h)$ & $\textrm{rate}_\tau$\\
	\hline
	& 32 & 2.2223E-1 & -- & 2.2227E-1 & -- & 2.2270E-1 & -- & 2.2270E-1 & -- \\
	& 64 & 1.0832E-1 & 1.0368 & 1.0832E-1 & 1.0371 & 1.0833E-1 & 1.0369 & 1.0833E-1 & 1.0369 \\
	DRSVD-LT & 128 & 5.4294E-2 & 0.9964 & 5.3160E-2 & 1.0269 & 5.3204E-2 & 1.0258 & 5.3206E-2 & 1.0257 \\
	& 256 & 3.1291E-2 & 0.7951 & 2.6318E-2 & 1.0143 & 2.6396E-2 & 1.0112 & 2.6400E-2 & 1.0112 \\
	& 512 & 2.4263E-2 & 0.3670 & 1.3119E-2 & 1.0044 & 1.3147E-2 & 1.0056 & 1.3154E-2 & 1.0051\\
	\hline
	& 32 & 2.2223E-1 & -- & 2.2227E-1 & -- & 2.2227E-1 & -- & 2.2227E-1 & --\\
	& 64 & 1.0832E-1 & 1.0368 & 1.0832E-1 & 1.0371 & 1.0833E-1 & 1.0369 &  1.0833E-1 & 1.0369\\
	DGN-LT & 128 & 5.4294E-2 & 0.9964 & 5.3160E-2 & 1.0269 & 5.3204E-2 & 1.0258 & 5.3206E-2 & 1.0257\\
	& 256 & 3.1291E-2 & 0.7951 & 2.6318E-2 & 1.0143 & 2.6396E-2 & 1.0112 & 2.6400E-2 & 1.0110\\
	& 512 & 2.4263E-2 & 0.3670 & 1.3119E-2 & 1.0044 & 1.3147E-2 & 1.0056 & 1.3154E-2 & 1.0051\\
	\hline
\end{tabular}
\label{tab7}
\end{center}
\end{table}

Tables \ref{tab7} and \ref{tab8} illustrate the numerical performance of solving Eq.~\eqref{eq4.2} using DRSVD-LT, DGN-LT, DRSVD-ST and DGN-ST methods. The experimental results indicate that for the Lie-Trotter splitting method (as shown in Table \ref{tab7}), the error performance of both DRSVD and DGN is comparable. As the rank-$r$ and the number of time steps $\hat{M}$ increase, the relative error gradually decreases, with the Lie-Trotter splitting method typically demonstrating the first-order convergence. In Table \ref{tab8}, we present the corresponding results for the Strang splitting method, which exhibits the anticipated second-order convergence for sufficiently high ranks.

\begin{table}[ht]\tabcolsep=3.0pt
\renewcommand{\arraystretch}{1.3} 
\begin{center}
\caption{Relative errors in Frobenius norm and temporal convergence orders of DRSVD-ST and DGN-ST ($\hat{N}=400$) for Eq.~\eqref{eq4.2}.}
\centering
\scriptsize
\begin{tabular}{cccccccccc}
	\hline
	& & \multicolumn{2}{c}{$r = 4$} & \multicolumn{2}{c}{$r = 6$} & \multicolumn{2}{c}{$r = 8$} 
	& \multicolumn{2}{c}{$r = 10$}\\
	\cmidrule(lr){3-4} \cmidrule(lr){5-6} \cmidrule(lr){7-8} \cmidrule(lr){9-10}
	Method & $\hat{M}$ & $\textrm{relerr}(\tau,h)$ & $\textrm{rate}_\tau$ 
	& $\textrm{relerr}(\tau,h)$ & $\textrm{rate}_\tau$
	& $\textrm{relerr}(\tau,h)$ & $\textrm{rate}_\tau$ & $\textrm{relerr}(\tau,h)$ & $\textrm{rate}_\tau$\\
	\hline
	& 32 & 8.8030E-2 & -- & 8.7378E-2 & -- & 8.7377E-2 & -- & 8.7377E-2 & -- \\
	& 64 & 3.0292E-2 & 1.5391 & 2.5967E-2 & 1.7506 & 2.5935E-2 & 1.7523 & 2.5935E-2 & 1.7523 \\
	DRSVD-ST & 128 & 1.9494E-2 & 0.6359 & 7.1414E-3 & 1.8624 & 6.8897E-3 & 1.9124 & 6.8890E-3 & 1.9126 \\
	& 256 & 2.0783E-2 & -0.0924 & 2.7755E-3 & 1.3635 & 1.7606E-3 & 1.9684 & 1.7530E-3 & 1.9745 \\
	& 512 & 2.2403E-2 & -0.1083 & 2.4010E-3 & 0.2091 & 4.7606E-4 & 1.8868 & 4.4038E-4 & 1.9930 \\
	\hline
	& 32 & 8.8030E-2 & -- & 8.7378E-2 & -- & 8.7377E-2 & -- & 8.7377E-2 & -- \\
	& 64 & 3.0292E-2 & 1.5391 & 2.5967E-2 & 1.7506 & 2.5935E-2 & 1.7523 & 2.5935E-2 & 1.7523\\
	DGN-ST & 128 & 1.9494E-2 & 0.6359 & 7.1414E-3 & 1.8624 & 6.8897E-3 & 	1.9124 & 6.8890E-3 & 1.9126\\
	& 256 & 2.0783E-2 & -0.0924 & 2.7755E-3 & 1.3635 & 1.7606E-3 & 	1.9684 & 1.7530E-3 & 1.9745\\
	& 512 & 2.2403E-2 & -0.1083 & 2.4010E-3 & 0.2091 & 4.7606E-4 & 	1.8868 & 4.4038E-4 & 1.9930\\
	\hline
\end{tabular}
\label{tab8}
\end{center}
\end{table}

\begin{figure}[H]	
\centering
\subfigure{
\includegraphics[width=2.08in,height=2.45in]{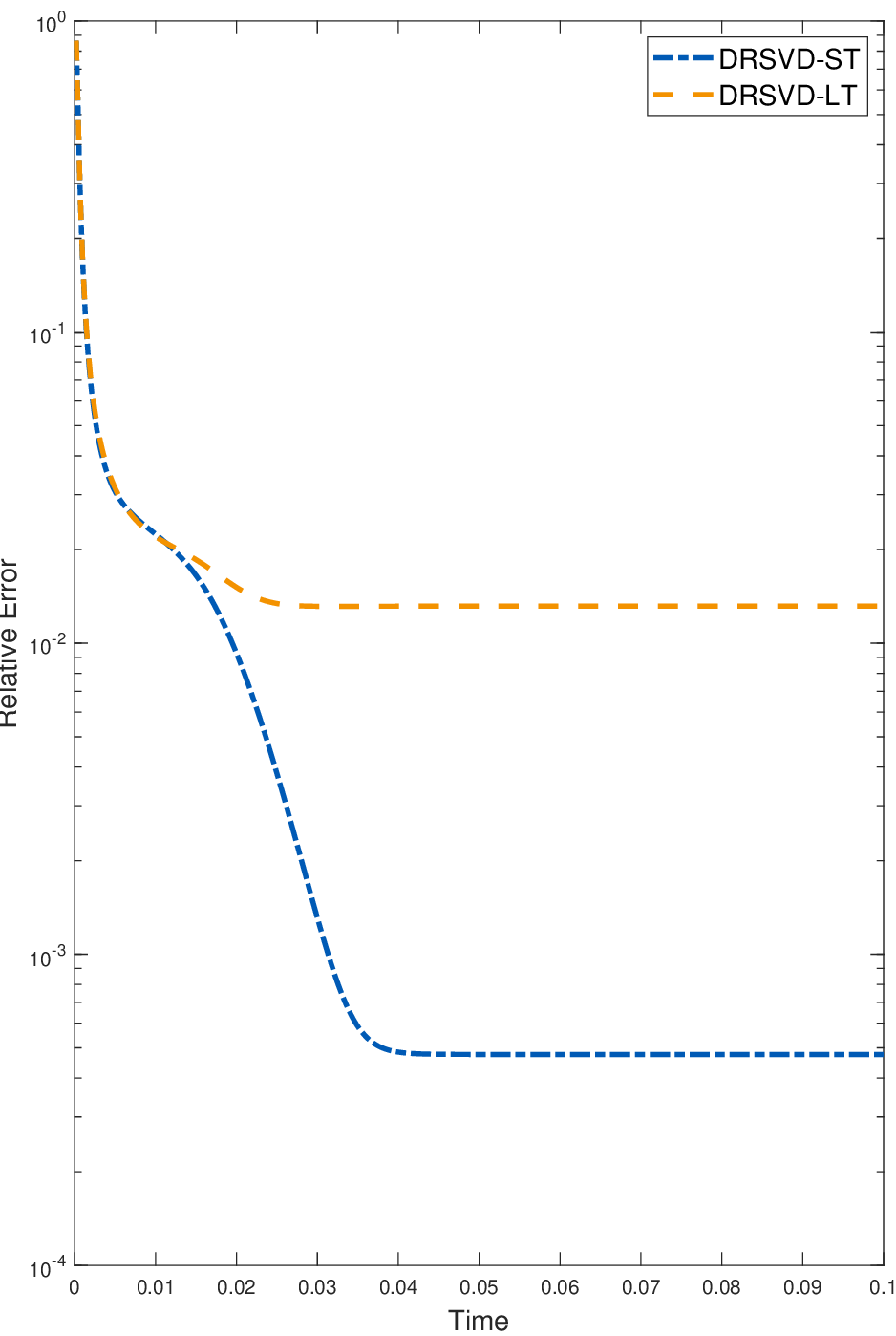}} \hspace{10mm}
\subfigure{
\includegraphics[width=2.08in,height=2.45in]{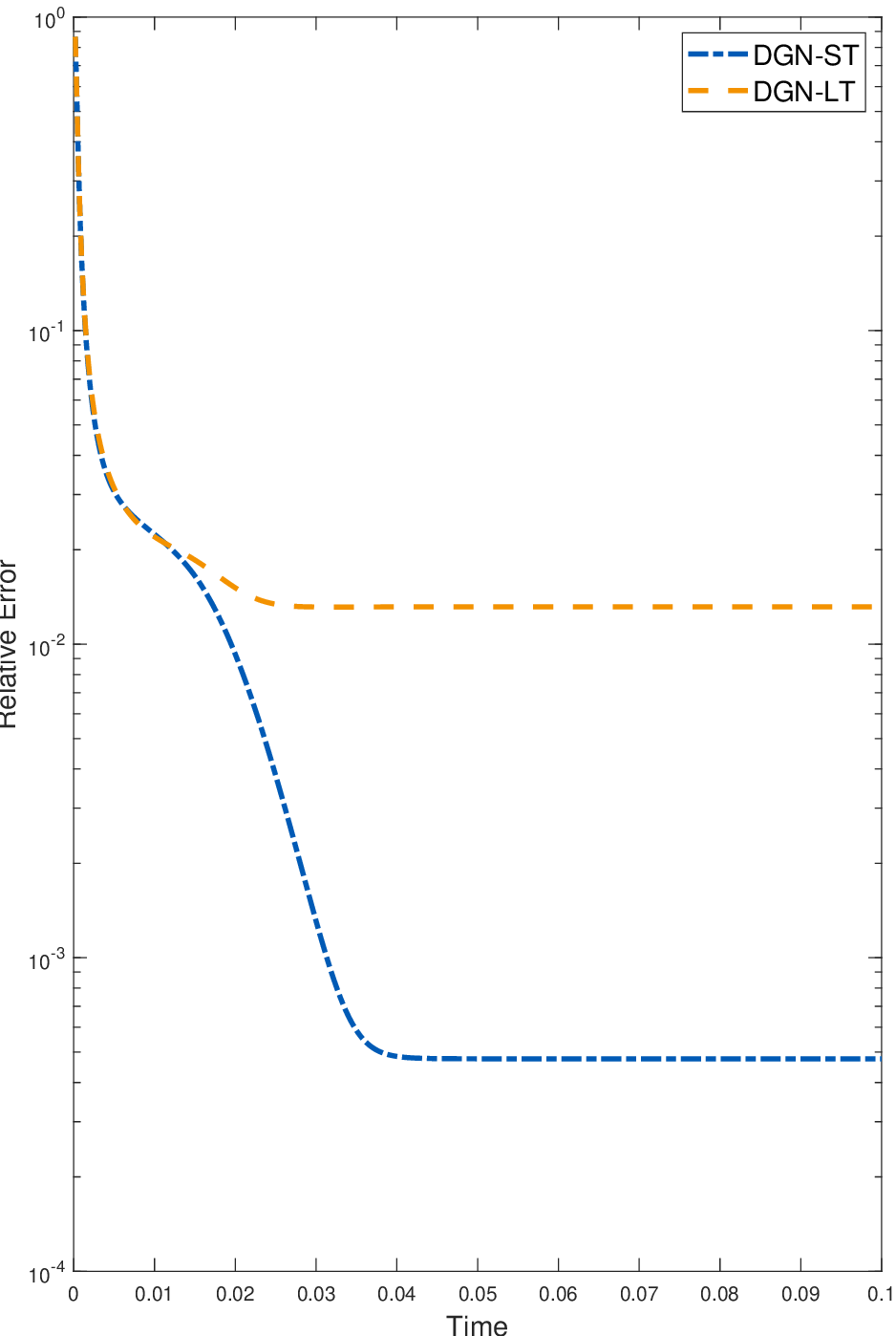}}			
\caption{Comparison of relative errors between Lie-Trotter splitting, Strang splitting and best rank approximation over time  for  Eq.~\eqref{eq4.2} for target rank $r = 8$ and the number of time steps $\hat{M} = 512$.}
\label{fig5}
\end{figure}

Fig.~\ref{fig5} illustrates the temporal evolution of errors for Lie-Trotter splitting and Strang splitting methods under the conditions of a fixed rank $r = 8$ and the number of time steps $\hat{M} = 512$. Specifically, the left subplot displays the results obtained using DRSVD to handle the nonlinear term, while the right subplot presents the results using DGN to address the nonlinear term. As anticipated, regardless of which nonlinear treatment (DRSVD or DGN) is employed, the error curve of the Strang splitting method consistently remains below that of the Lie-Trotter splitting method, confirming its inherent advantage in achieving the second-order accuracy.

\begin{table}[ht]\tabcolsep=3.0pt
\renewcommand{\arraystretch}{1.3} 
\begin{center}
\caption{Relative errors in Frobenius norm and temporal convergence orders of ADRSVD and ADGN ($\hat{N}=400$) for Eq.~\eqref{eq4.2}.}
\centering
\scriptsize
\begin{tabular}{ccccccccc}
	\hline
	& \multicolumn{2}{c}{ADRSVD-LT} & \multicolumn{2}{c}{ADGN-LT} & \multicolumn{2}{c}{ADRSVD-ST} & \multicolumn{2}{c}{ADGN-ST}\\
	\cmidrule(lr){2-3}\cmidrule(lr){4-5}\cmidrule(lr){6-7}\cmidrule(lr){8-9}
	$\hat{M}$ & $\textrm{relerr}(\tau,h)$ & $\textrm{rate}_\tau$ 
	& $\textrm{relerr}(\tau,h)$ & $\textrm{rate}_\tau$& $\textrm{relerr}(\tau,h)$ & $\textrm{rate}_\tau$ 
	& $\textrm{relerr}(\tau,h)$ & $\textrm{rate}_\tau$\\
	\hline
	32 & 2.2227E-1& -- & 2.2227E-1 & -- & 8.7377E-2 & -- & 8.7377E-2 & -- \\
	64 & 1.0833E-1 & 1.0369 & 1.0833E-1 & 1.0369 & 2.5935E-2 & 1.7523 & 2.5935E-2 & 1.7523\\ 
	128 & 5.3206E-2 & 1.0257 & 5.3206E-2 & 1.0257 & 6.8889E-3 & 1.9126 & 6.8889E-3 & 1.9126\\
	256 & 2.6400E-2 & 1.0110 & 2.6400E-2 & 1.0110 & 1.7529E-3 & 1.9745 & 1.7529E-3 & 1.9745\\
	512 & 1.3154E-2 & 1.0051 & 1.3154E-2 & 1.0051 & 4.4026E-4 & 1.9934 & 4.4026E-4 & 1.9934\\
	\hline
\end{tabular}
\label{tab9}
\end{center}
\end{table}

\begin{figure}[H]
\centering
\subfigure
{\includegraphics[width=2.1in,height=2.1in]{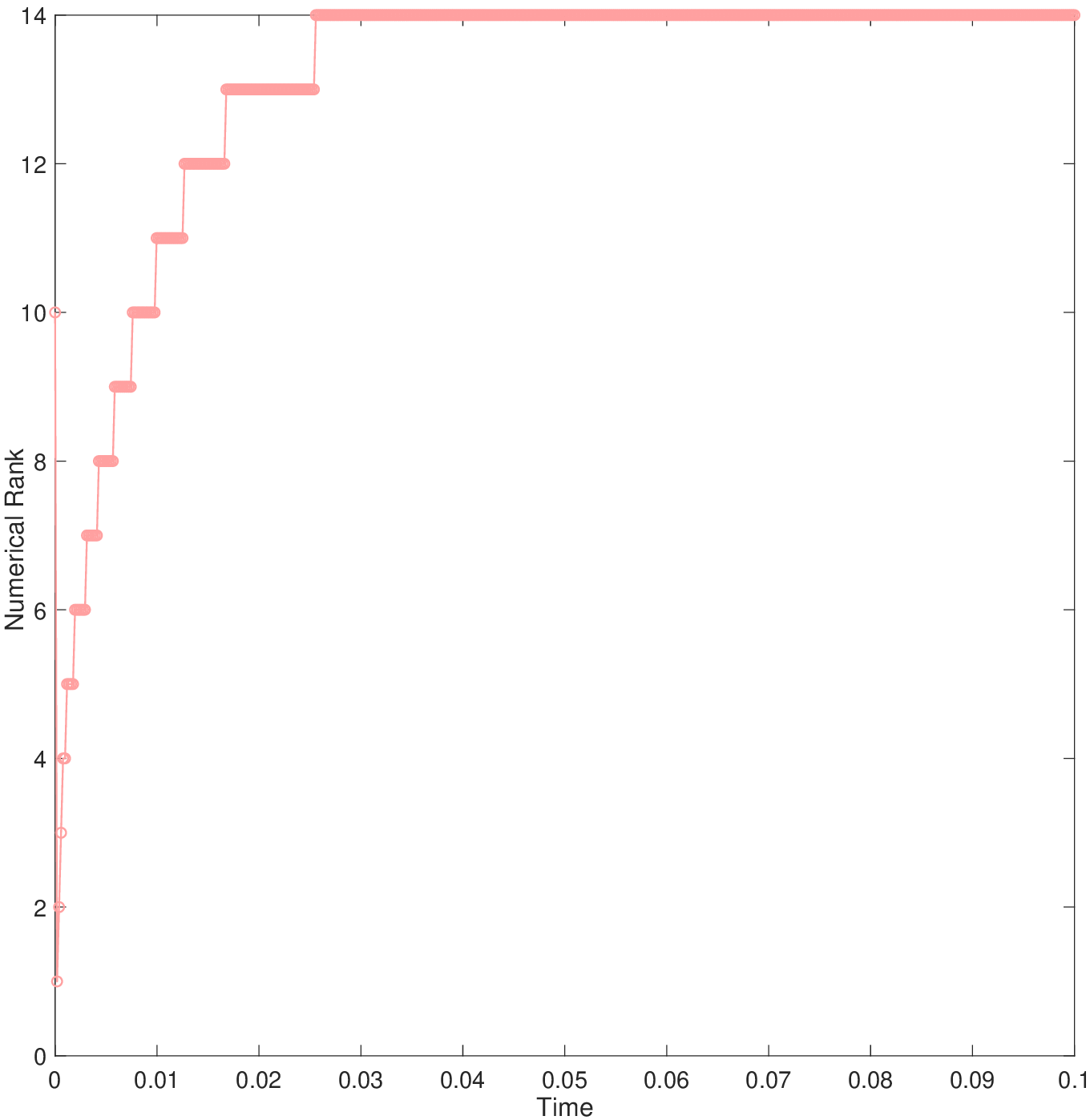}} \hspace{2mm}
\subfigure
{\includegraphics[width=2.1in,height=2.1in]{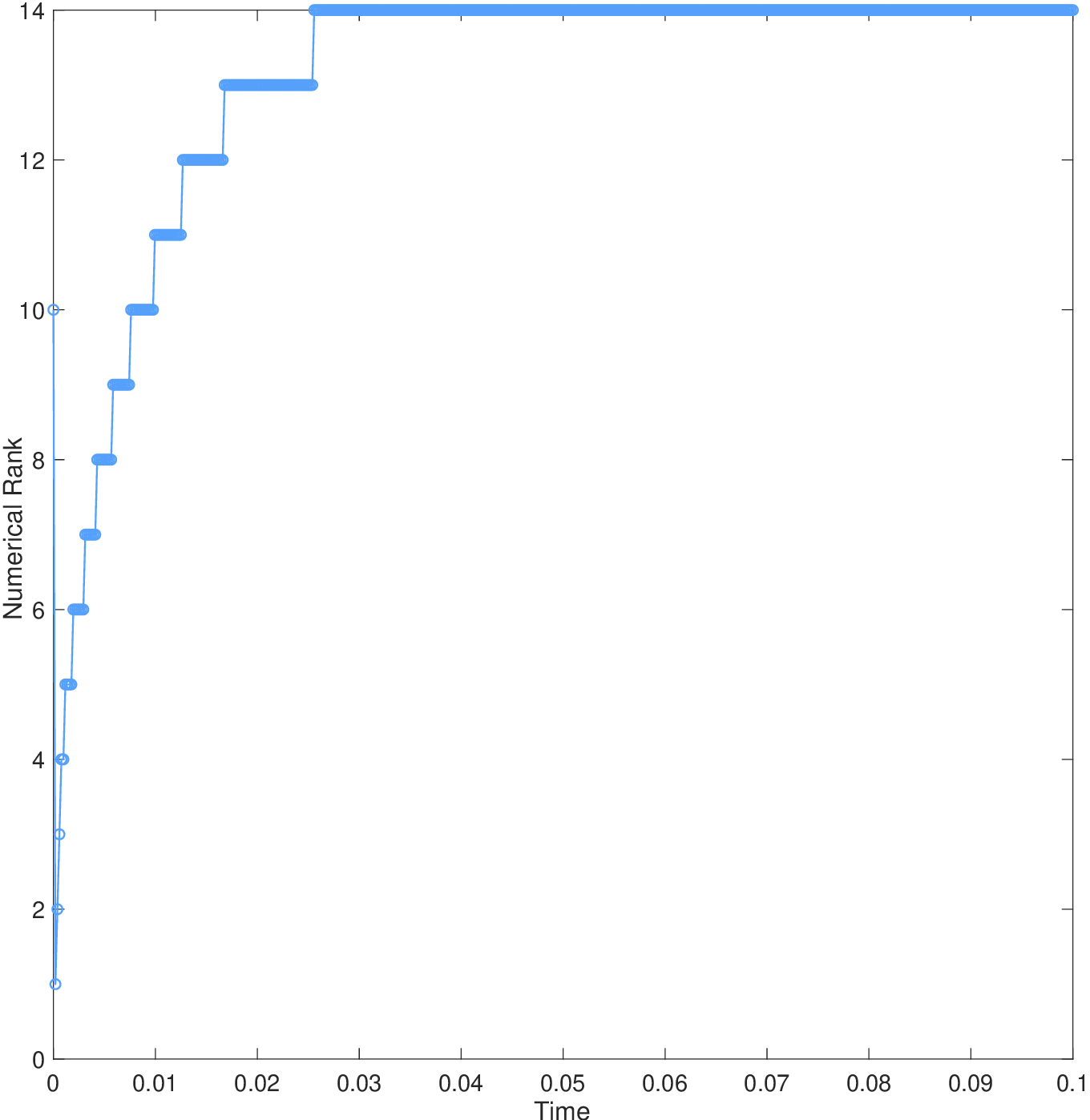}}
\caption{ Numerical convergence order evolution of Eq.~\eqref{eq4.2} for splitting methods combined with low-rank solvers $(\hat{M},\hat{N})=(512,400)$.
Left: ADRSVD-LT; 
Right: ADRSVD-ST.}
\label{fig6}
\end{figure}

Table \ref{tab9} presents the numerical results obtained using ADRSVD-LT, ADGN-LT, ADRSVD-ST and ADGN-ST methods. For the truncation step, the relative tolerance is set to \(\varphi = 10^{-8}\) and the absolute tolerance to \(10^{-12}\). 
The rangefinder step employs a tolerance of \(10^{-8}\). The datas in tables clearly indicate that, under the same splitting method, the relative errors (relerr) and convergence orders (rate) of ADRSVD and ADGN are nearly identical. This suggests that the differences in these low-rank methods do not significantly affect the results. On the other hand, the relative errors decrease steadily with increasing the number of time steps $\hat{M}$ and strang splitting performs exceptionally well with fine time steps; its convergence order approaches the theoretical second-order, and the relative error is low, demonstrating a significant advantage over the Lie-Trotter splitting method. Furthermore, in Fig. \ref{fig6}, we show the change in rank over time with a fixed number of time steps $\hat{M} = 512$.

\subsection{Simulation}\label{sec4.2}
In this section, to validate the numerical stability of the ARSVD-LT, ADGN-LT, ADRSVD-ST and ADGN-ST methods in long-term simulations, we focus on the following Allen-Cahn equation with logarithmic Flory--Huggins potential \cite{wu2023second}:
\begin{equation}\label{eq4.4}
\begin{cases}
\partial_tu=\varepsilon^2\Delta u-f(u),\quad(t,\mathbf{x})\in(0,\infty)\times\Omega \\
F(u)=\theta\left[u\mathrm{ln}u+(1-u)\mathrm{ln}(1-u)\right]+2\theta_cu(1-u),\mathrm{~for}\quad0<u<1 \\
\left.u\right|_{t=0}=u^0, & 
\end{cases}
\end{equation}
where $u$ represents a real-valued function, ${\mathbf{\varepsilon}}$ is defined as the mobility coefficient, \(F(u)\) signifies the logarithmic Flory-Huggins potential, with \(f(u) = F^\prime(u)\). $\theta$ and $\theta_{c}$ are the absolute temperature and critical temperature respectively. Here we set ${\mathbf{\varepsilon}} = 0.1$, $\theta = 0.8$ and $\theta_c = 1$.

Besides, we select two different initial values for featuring high-curvature interfaces:

(1)~\textbf{Star-shaped:} This is borrowed from \cite{mu2016new}, where 

\begin{equation*}
u(x,y,0) = \tanh \left( \frac{0.25 + 0.1 \cos(6 \tilde{\theta}) -
\sqrt{(x-.5)^2 + (y - 0.5^2)}}{\varepsilon \sqrt{2}} \right),
\end{equation*}
with $(x,y) \in [0,1]^2$ and
\begin{equation*}
\tilde{\theta} =
\begin{cases}
\tan^{-1} \left( \frac{y - 0.5}{x - 0.5} \right), & x > 0.5, \\
\pi + \tan^{-1} \left( \frac{y - 0.5}{x - 0.5} \right), & \mathrm{others}.
\end{cases}	
\end{equation*}

(2)~\textbf{Butterﬂy-shaped:} This employs the high-curvature, butterfly-shaped interface construction initial value proposed in \cite{mu2016new}, which is parameterized by the following equation:
\begin{equation*}
\begin{cases}
x(\hat{\theta})=k(a+b\cos(c\hat{\theta})\sin(d\hat{\theta}))\cos(\hat{\theta}), \\
y(\hat{\theta})=k(a+b\cos(c\hat{\theta})\sin(d\hat{\theta}))\sin(\hat{\theta}), & 
\end{cases}
\end{equation*}
where $\hat{\theta}\in[0,2\pi]$. In this case, we set $k = 1.8$, $a = 0.40178$, $b = 0.30178$, $c = 6.0$, and $d = 3.0$.

\begin{figure}[H]
\setlength{\tabcolsep}{0.2pt}
\centering
\begin{tabular}{m{0.5cm}<{\centering} m{3cm}<{\centering} m{3cm}<{\centering} m{3cm}<{\centering} m{3cm}<{\centering}}
& $T = 0$ & $T = 10$ & $T = 50$ & $T = 100$ \\
\rotatebox{90}{ADRSVD-LT} &
\includegraphics[width=1.2in,height=1.2in]{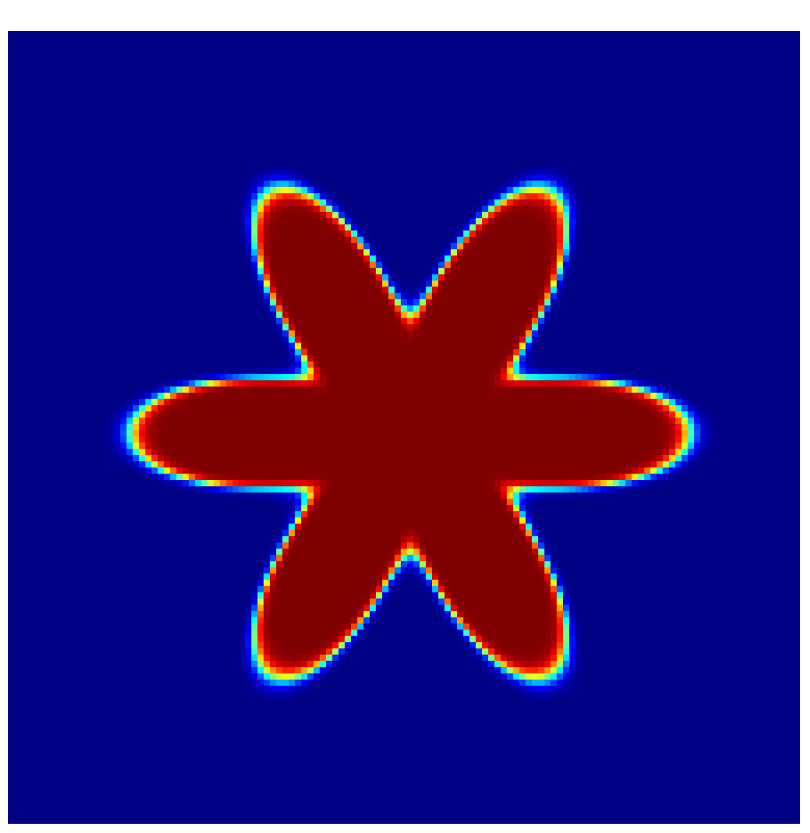} &
\includegraphics[width=1.2in,height=1.2in]{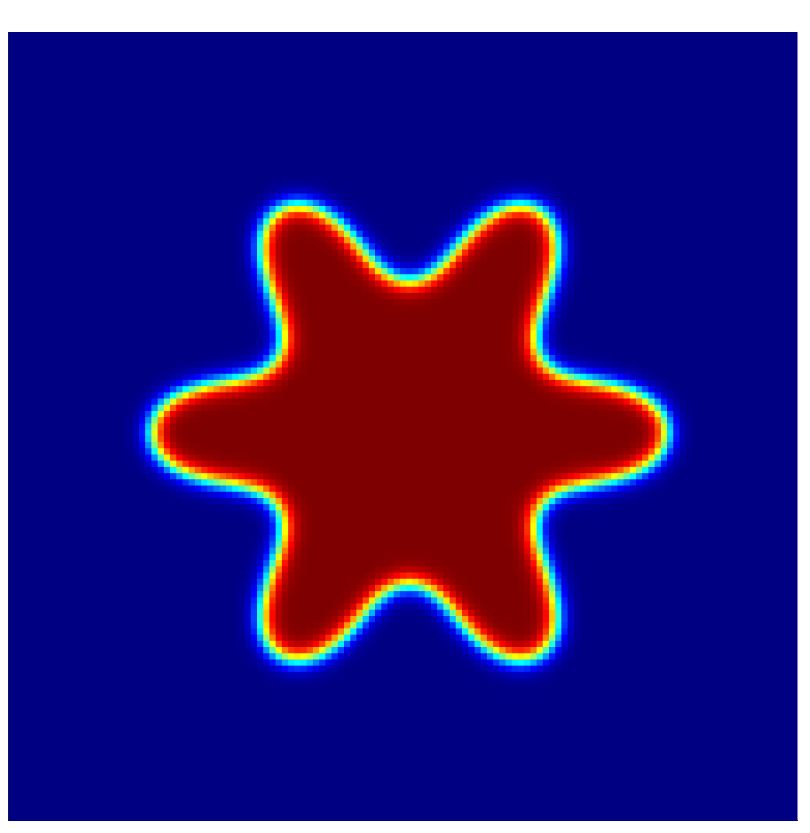} &
\includegraphics[width=1.2in,height=1.2in]{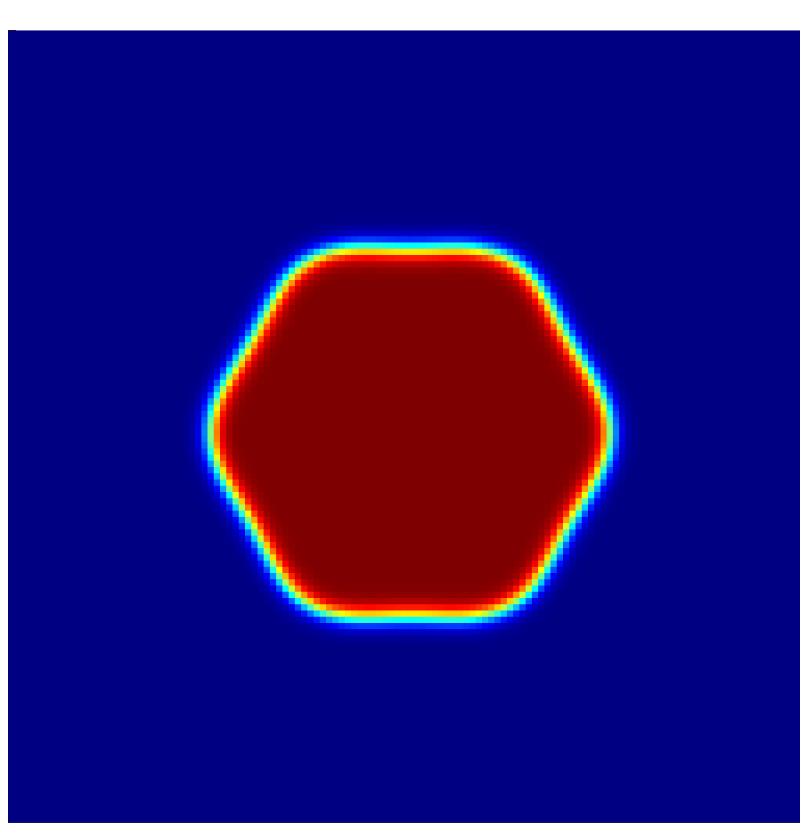} &
\includegraphics[width=1.2in,height=1.2in]{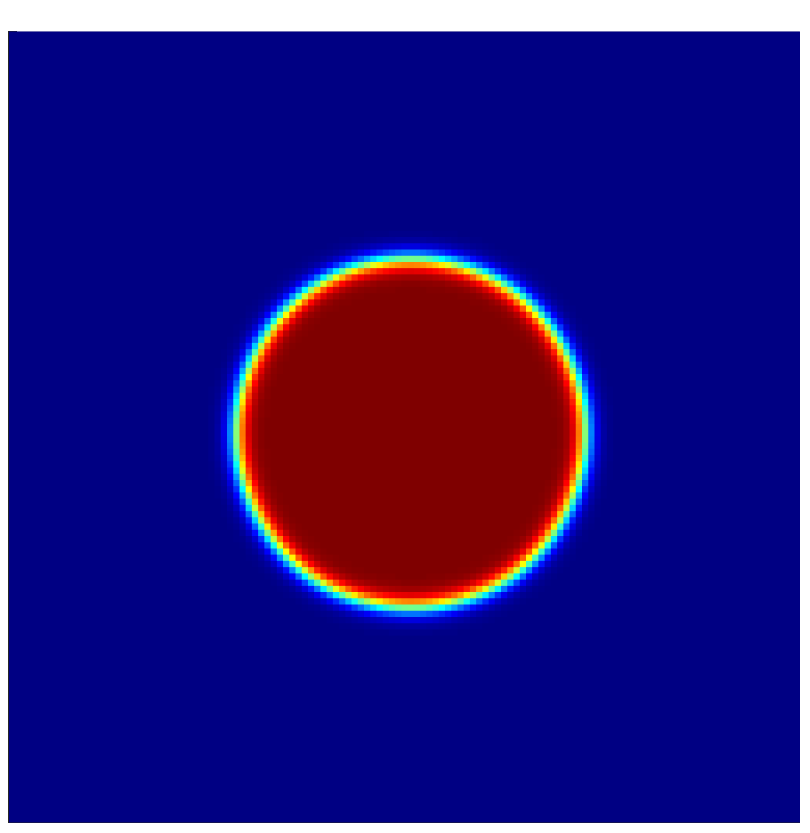} \\
\rotatebox{90}{ADGN-LT} &
\includegraphics[width=1.2in,height=1.2in]{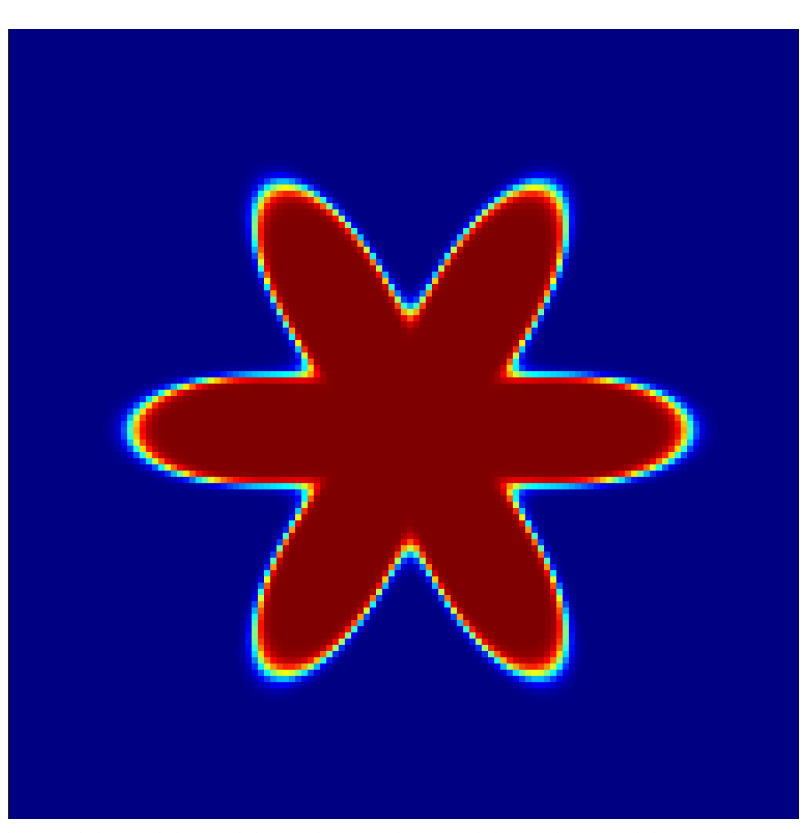} &
\includegraphics[width=1.2in,height=1.2in]{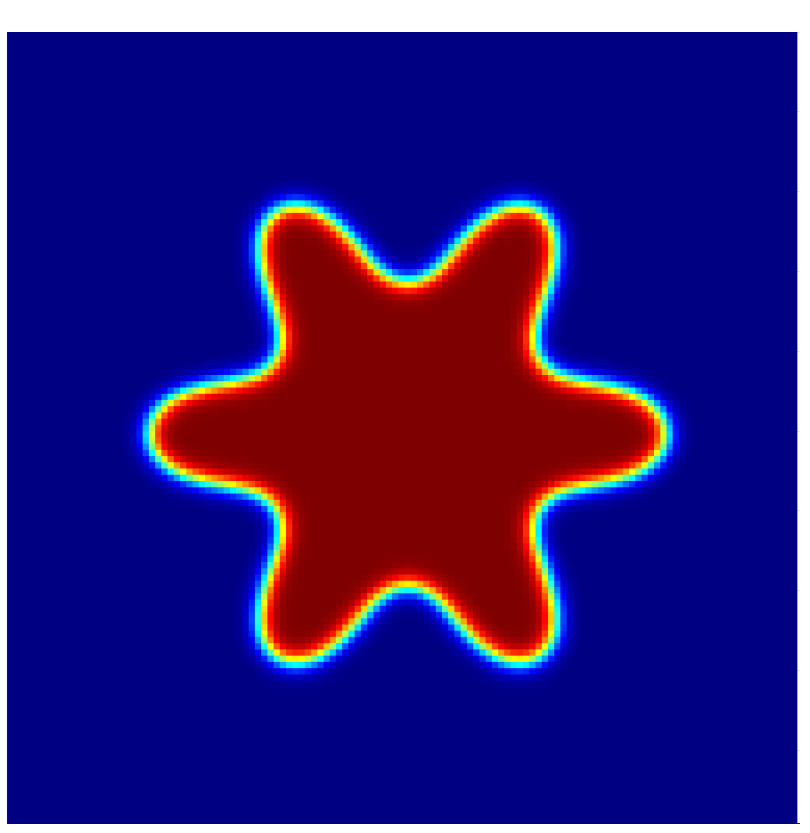} &
\includegraphics[width=1.2in,height=1.2in]{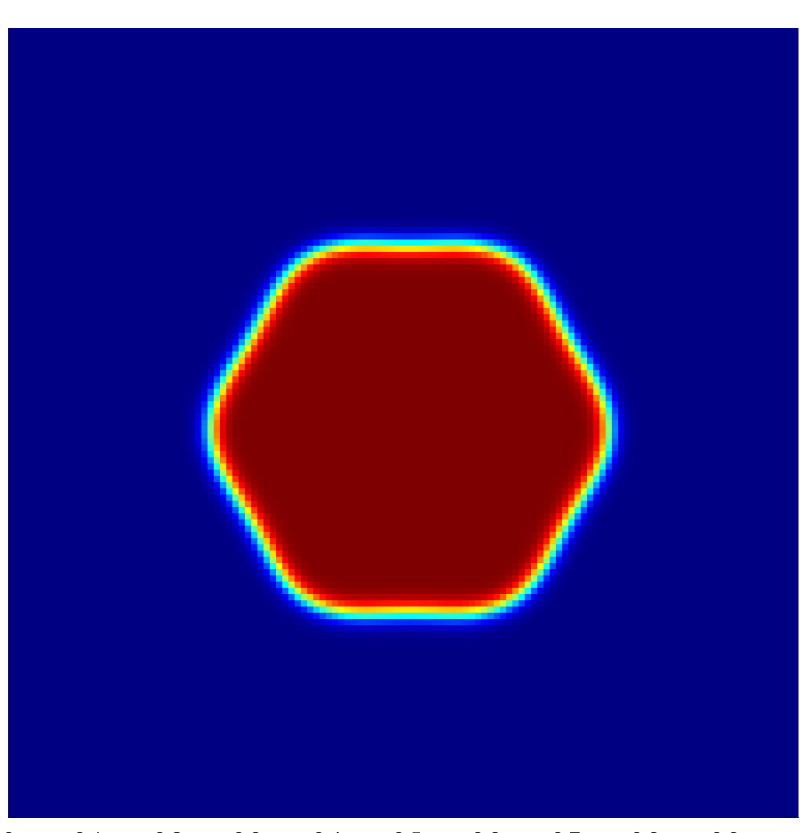} &
\includegraphics[width=1.2in,height=1.2in]{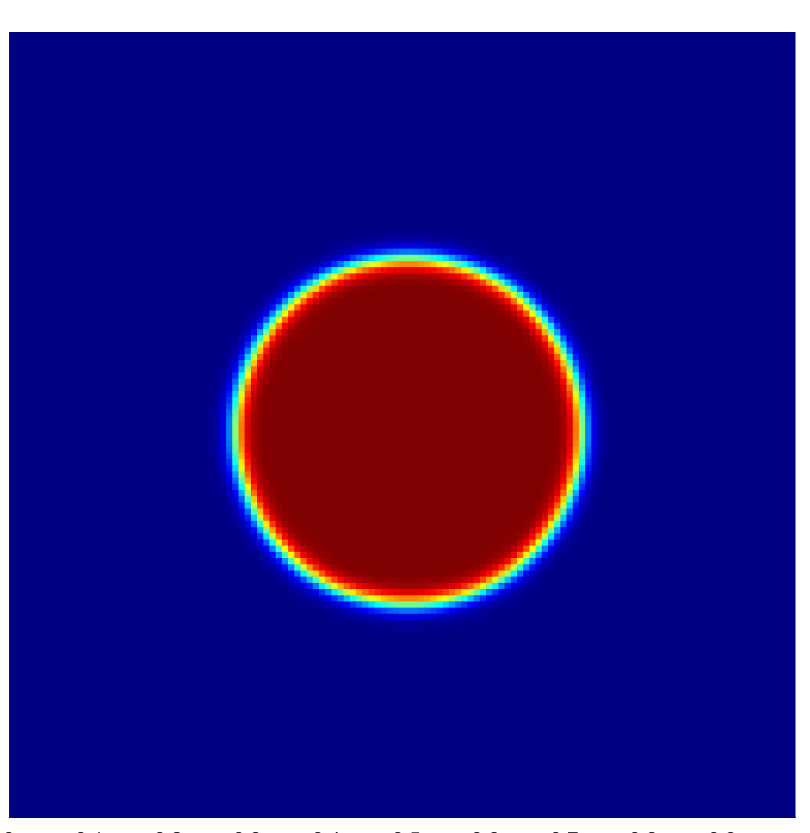} \\
\rotatebox{90}{ADRSVD-ST} &
\includegraphics[width=1.2in,height=1.2in]{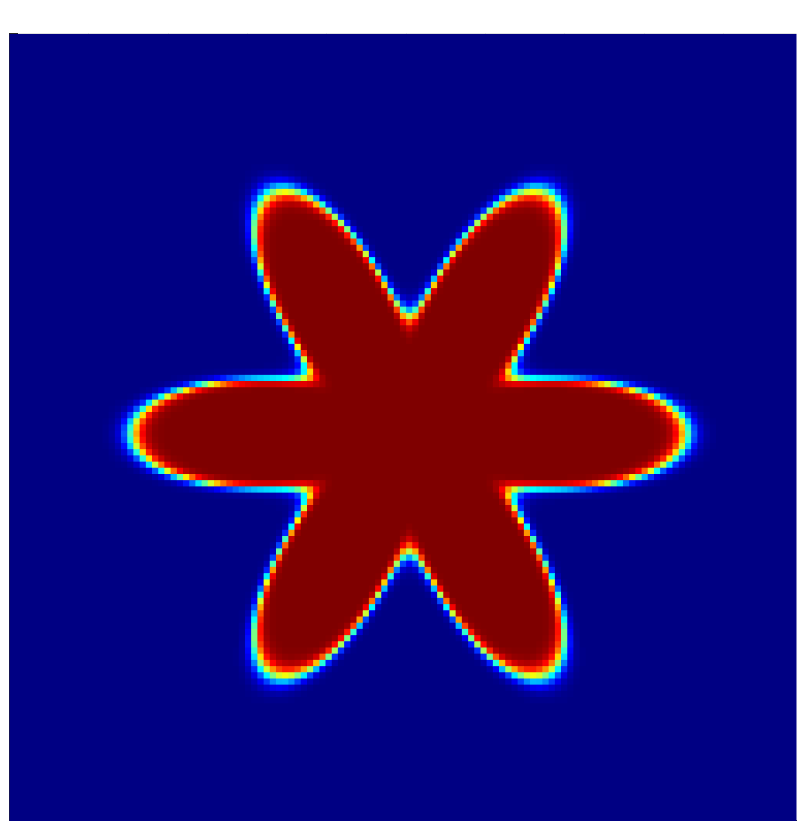} &
\includegraphics[width=1.2in,height=1.2in]{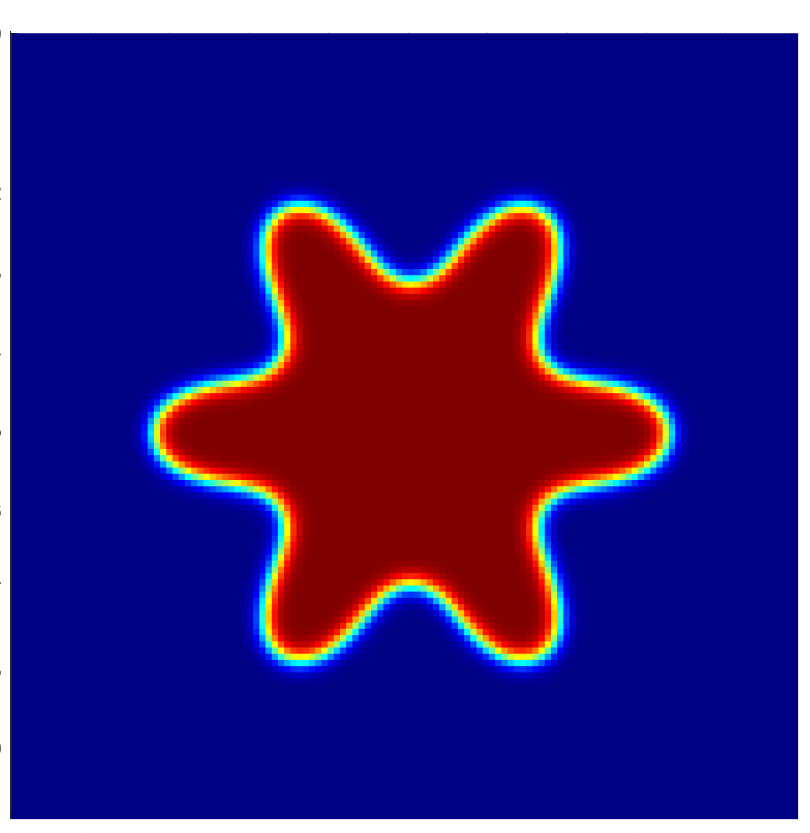} &
\includegraphics[width=1.2in,height=1.2in]{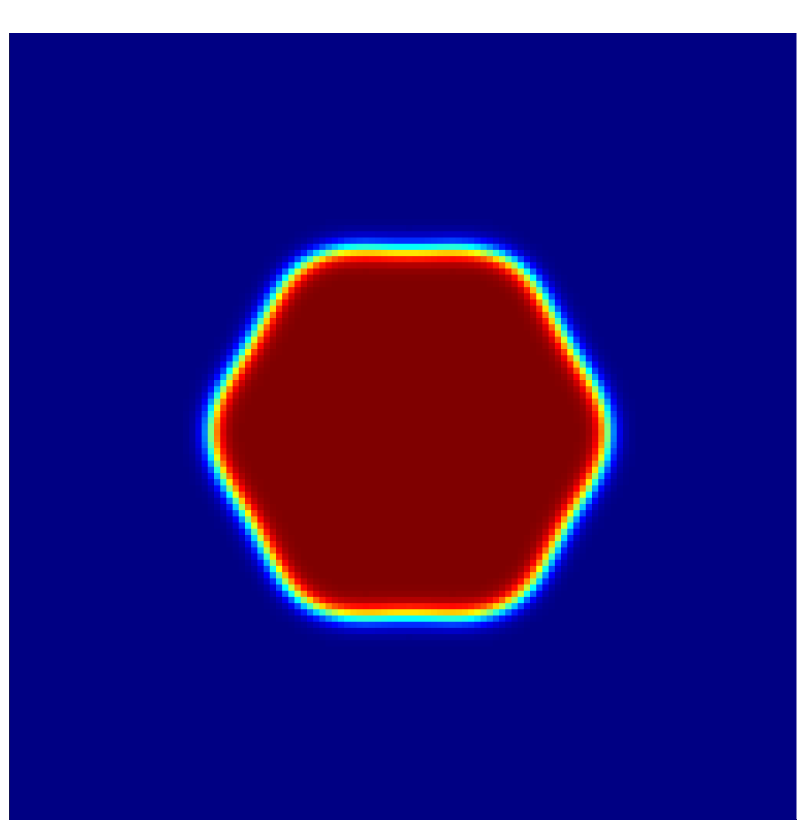} &
\includegraphics[width=1.2in,height=1.2in]{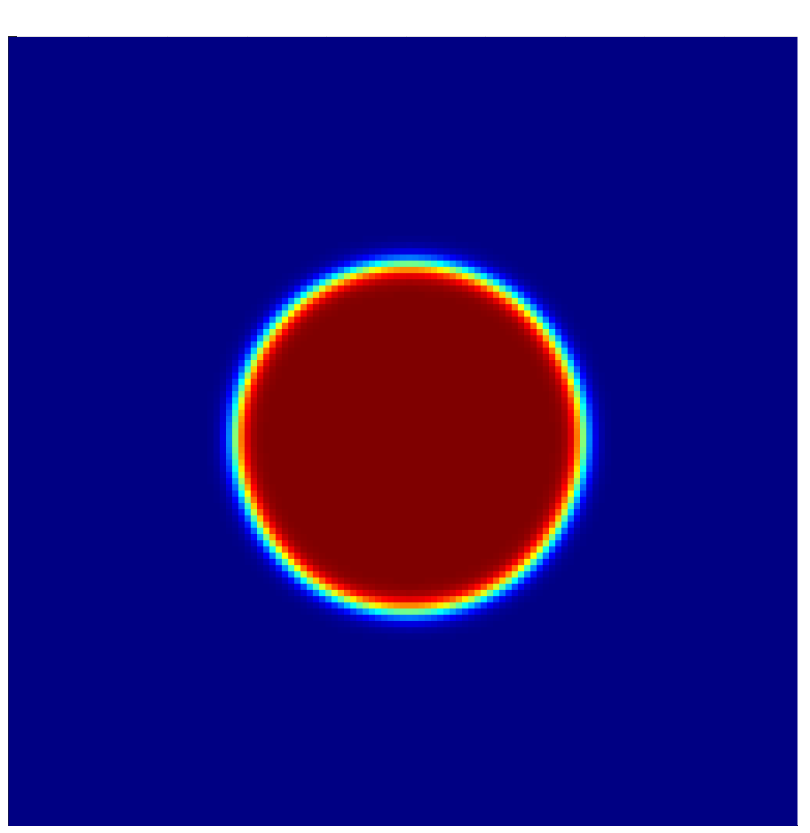}\\
\rotatebox{90}{ADGN-ST} &
\includegraphics[width=1.2in,height=1.2in]{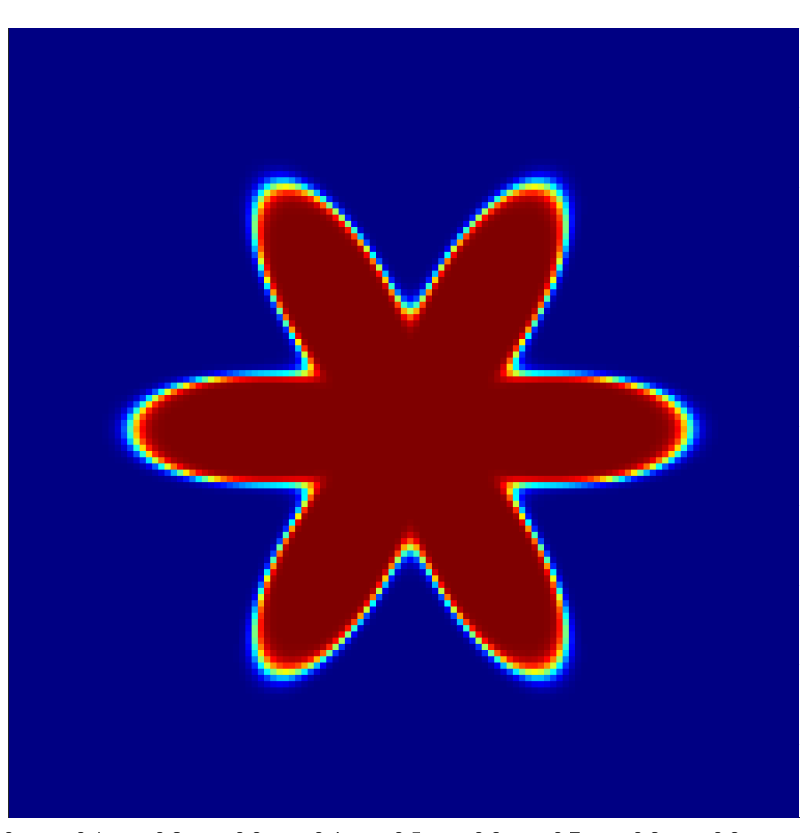} &
\includegraphics[width=1.2in,height=1.2in]{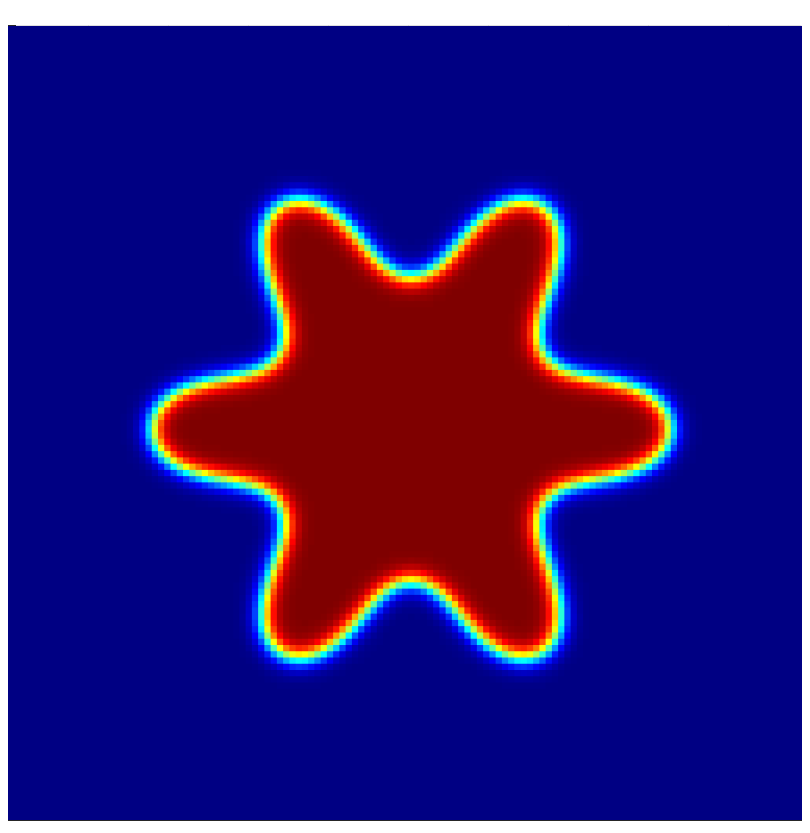} &
\includegraphics[width=1.2in,height=1.2in]{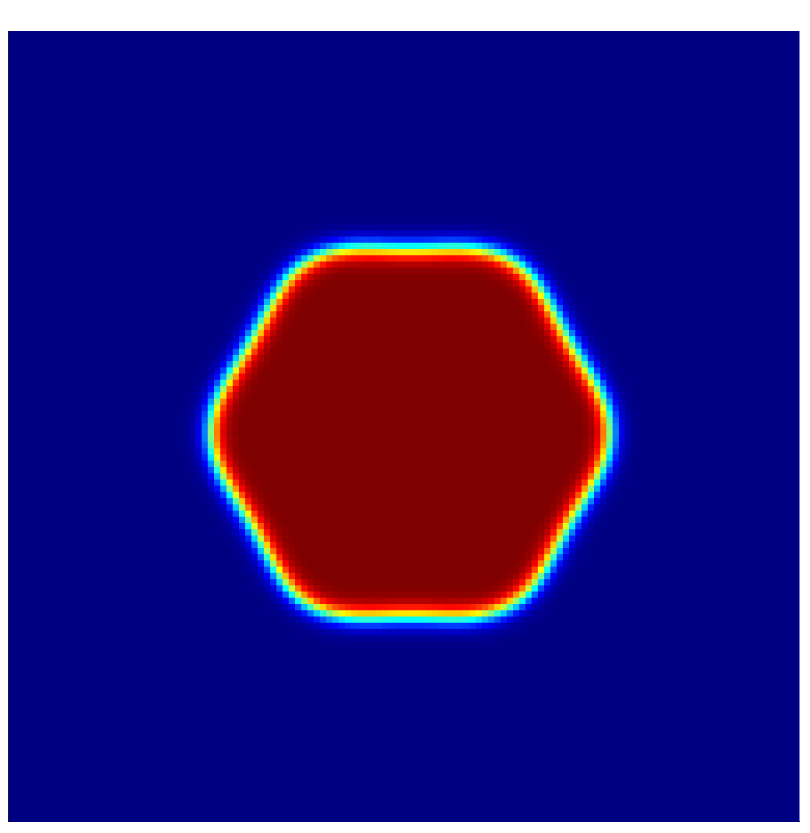} &
\includegraphics[width=1.2in,height=1.2in]{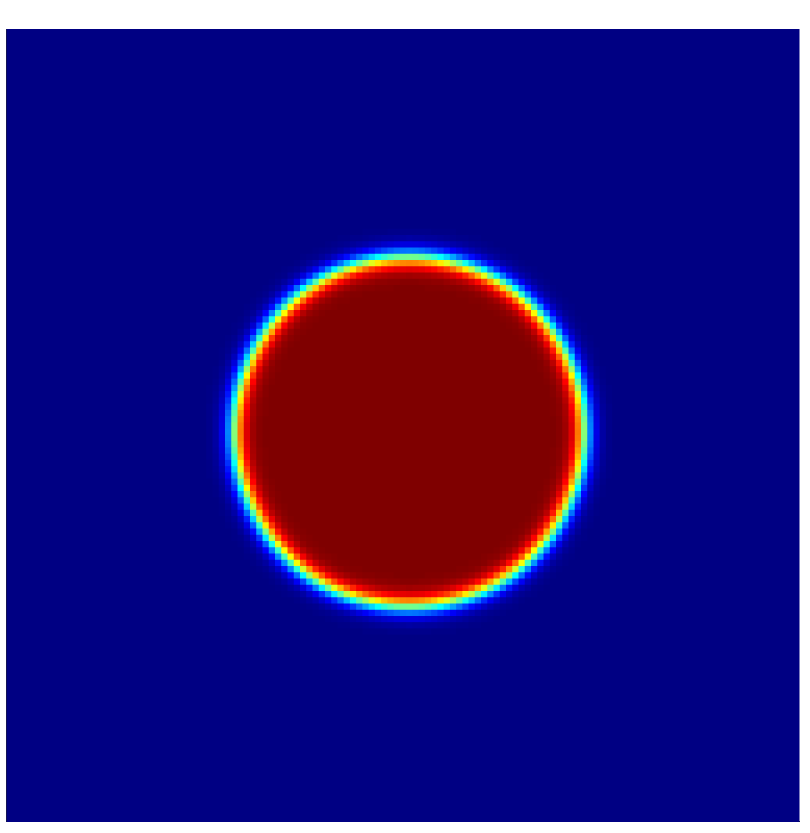} \\
\end{tabular}
\caption{Comparisons of ADRSVD-LT, ADGN-LT, ADRSVD-ST and ADGN-ST: dynamic evolution of star-shaped solutions for Eq.~\eqref{eq4.4} with $\hat{M} = 128$.
Top two rows: Lie-Trotter splitting; Bottom two rows: Strang splitting.}
\label{fig7}
\end{figure}

\begin{figure}[H]
\setlength{\tabcolsep}{0.2pt}
\centering
\begin{tabular}{m{0.5cm}<{\centering} m{3cm}<{\centering} m{3cm}<{\centering} m{3cm}<{\centering} m{3cm}<{\centering}}
& $T = 0$ & $T = 1$ & $T = 5$ & $T = 10$ \\
\rotatebox{90}{ADRSVD-LT} &
\includegraphics[width=1.2in,height=1.2in]{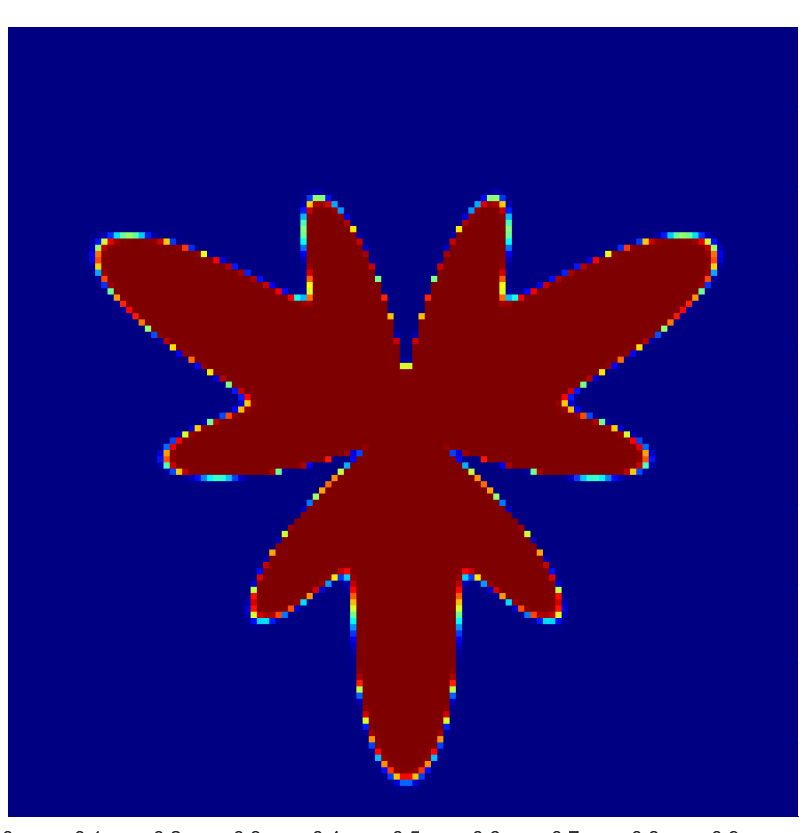} &
\includegraphics[width=1.2in,height=1.2in]{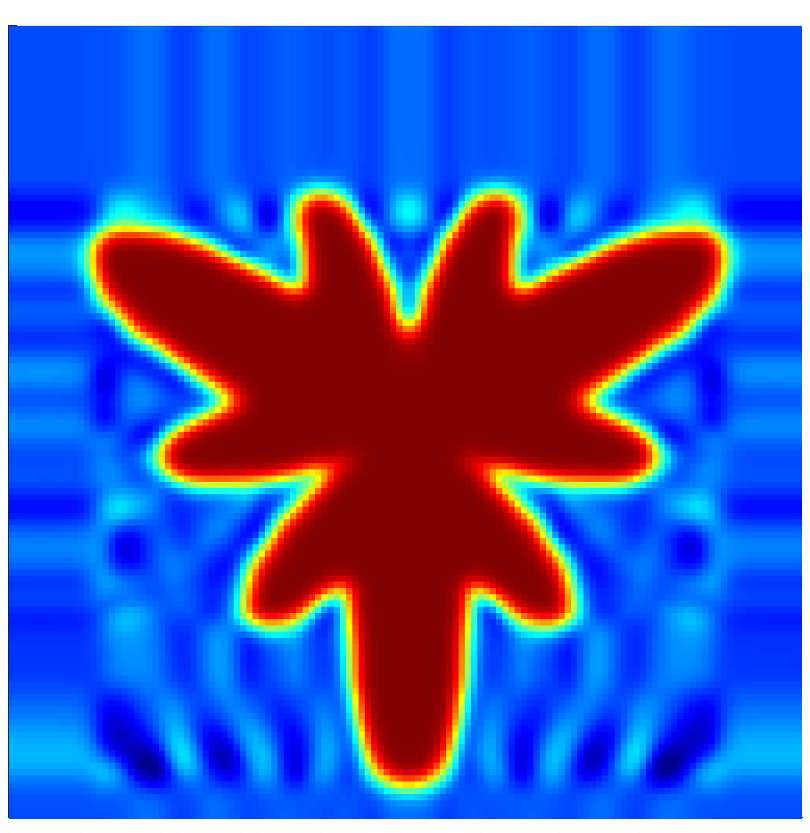} &
\includegraphics[width=1.2in,height=1.2in]{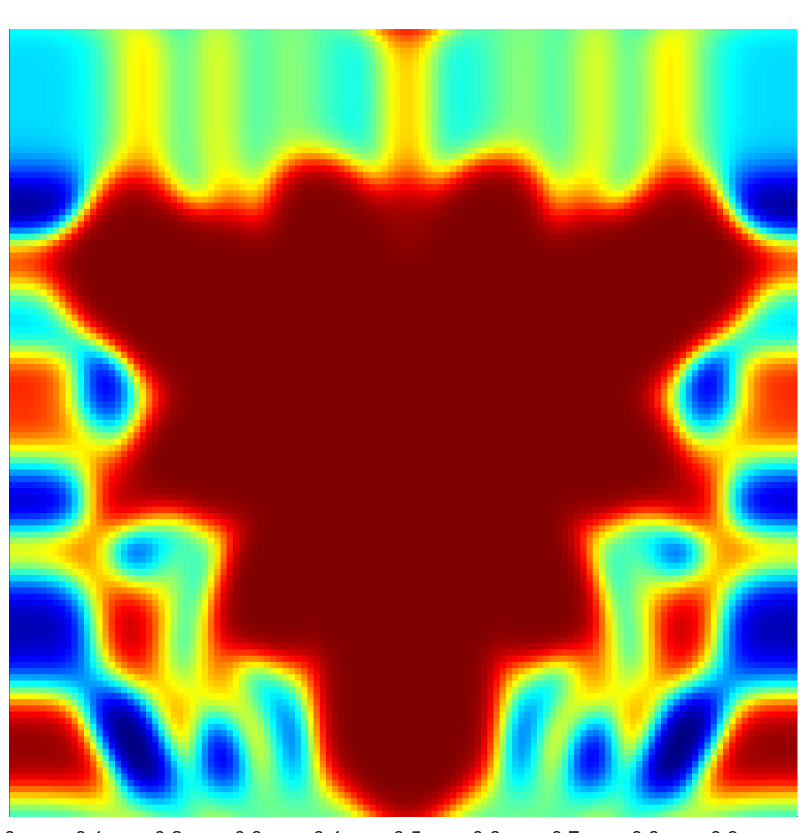} &
\includegraphics[width=1.2in,height=1.2in]{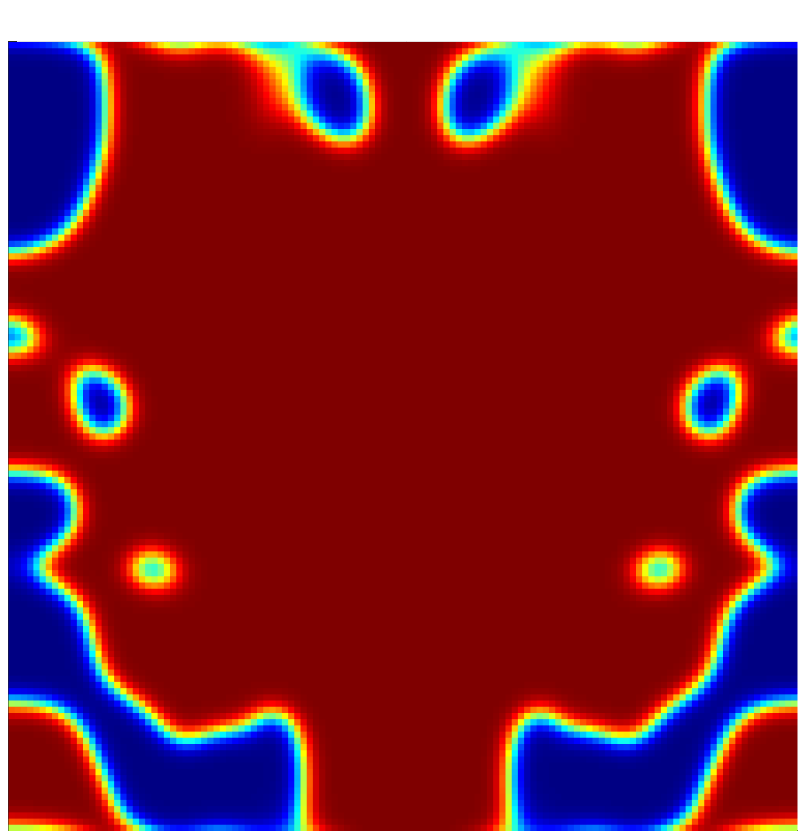} \\
\rotatebox{90}{ADGN-LT} &
\includegraphics[width=1.2in,height=1.2in]{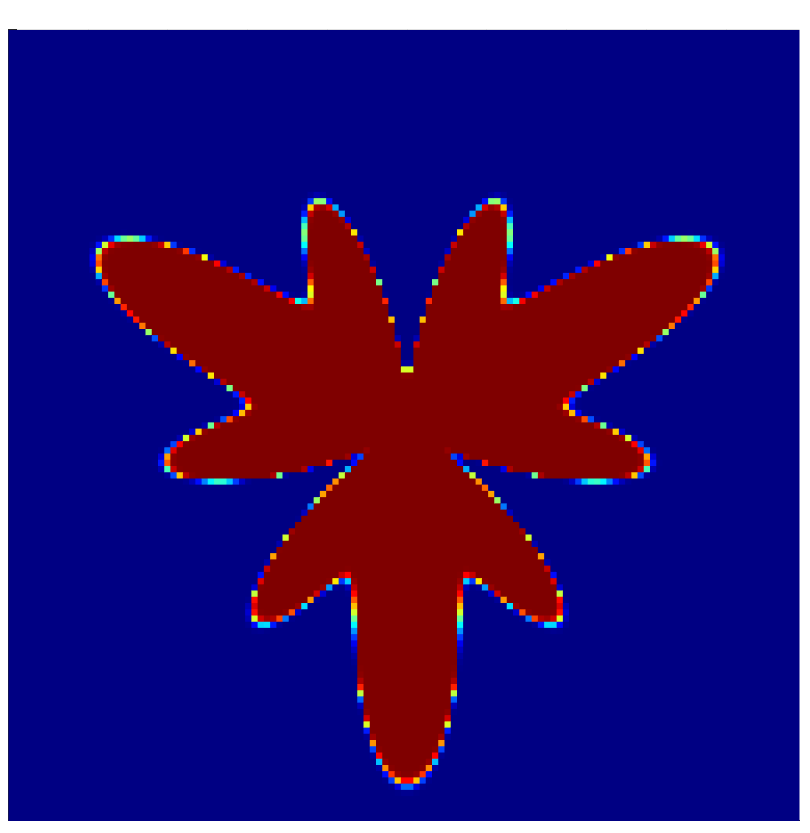} &
\includegraphics[width=1.2in,height=1.2in]{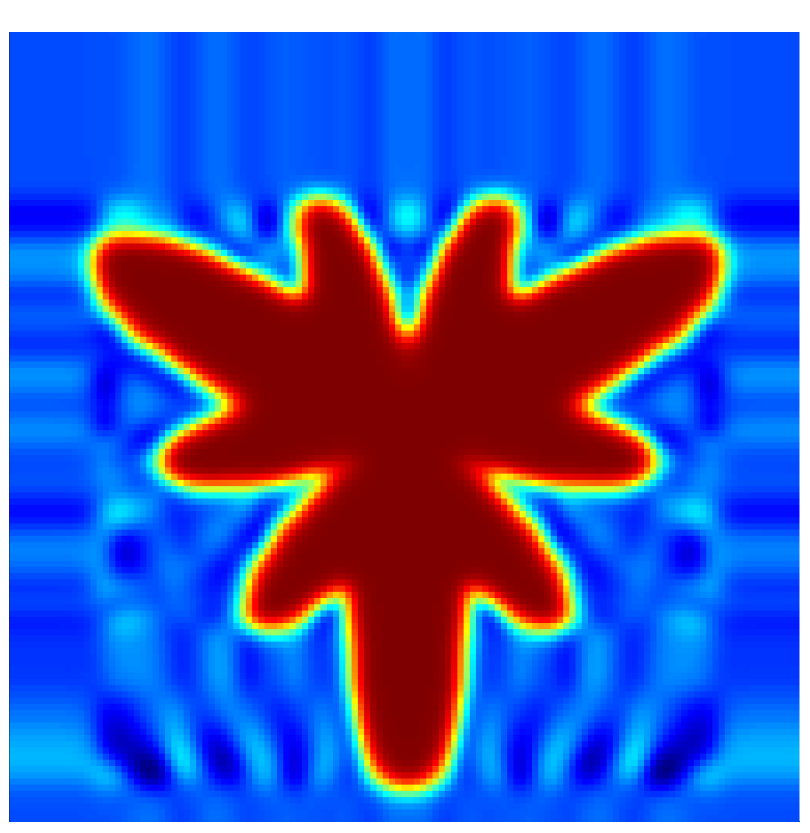} &
\includegraphics[width=1.2in,height=1.2in]{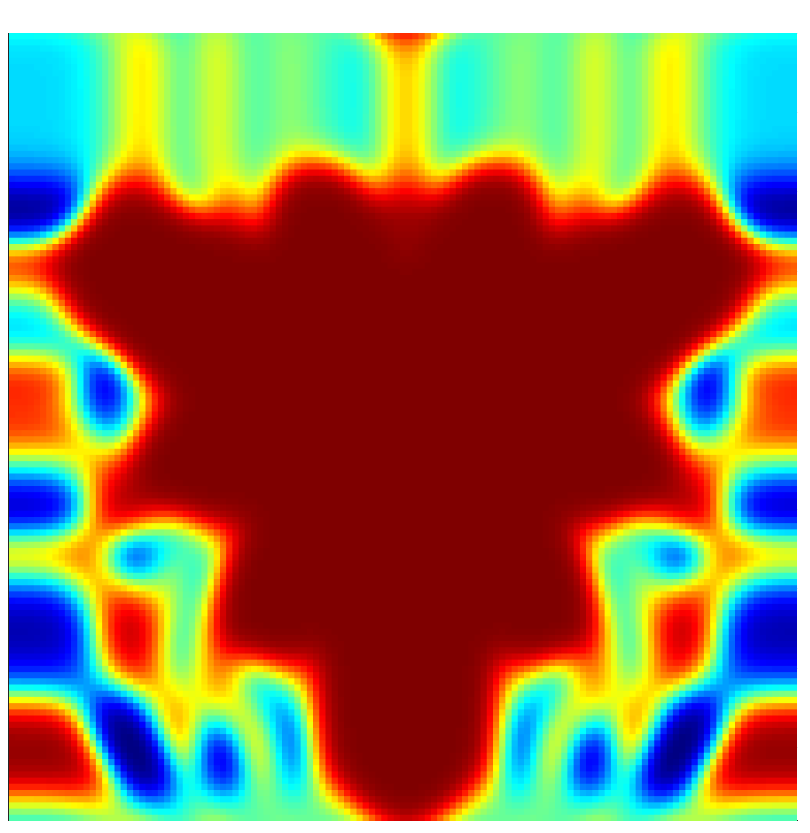} &
\includegraphics[width=1.2in,height=1.2in]{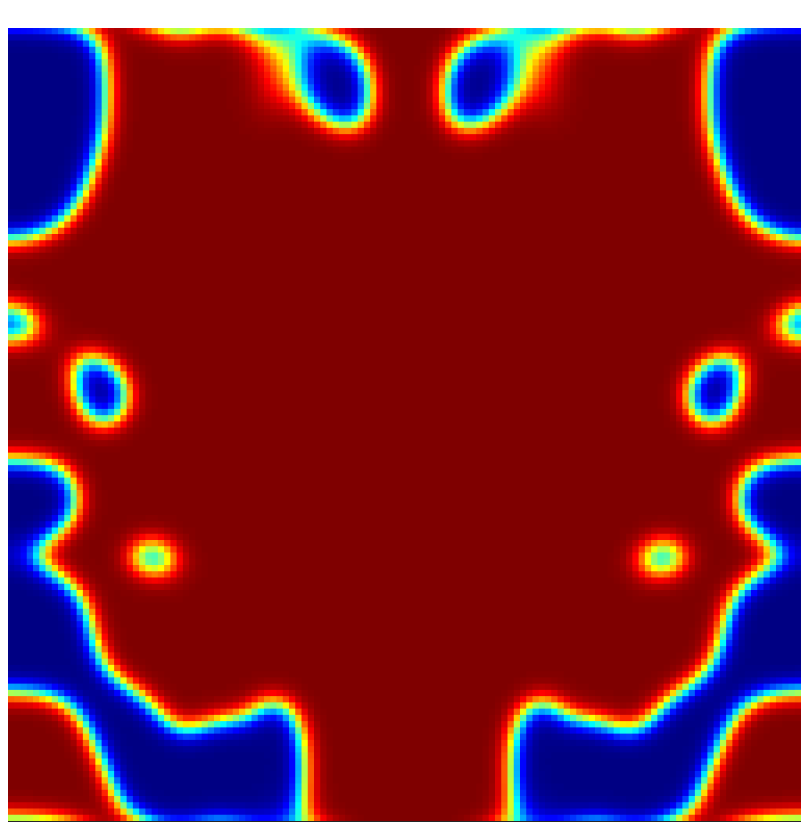} \\
\rotatebox{90}{ADRSVD-ST} &
\includegraphics[width=1.2in,height=1.2in]{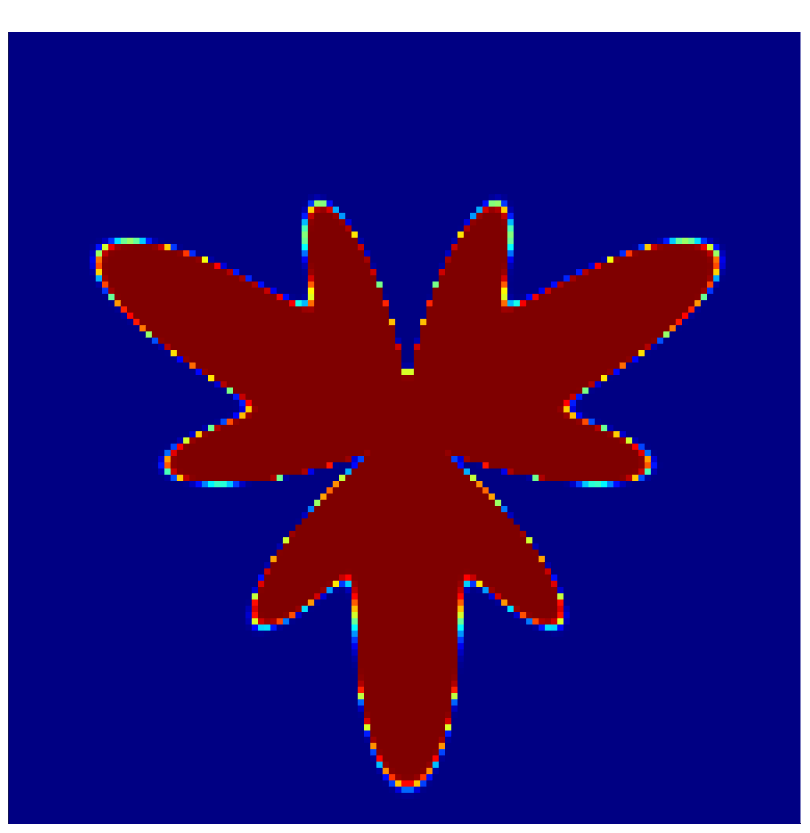} &
\includegraphics[width=1.2in,height=1.2in]{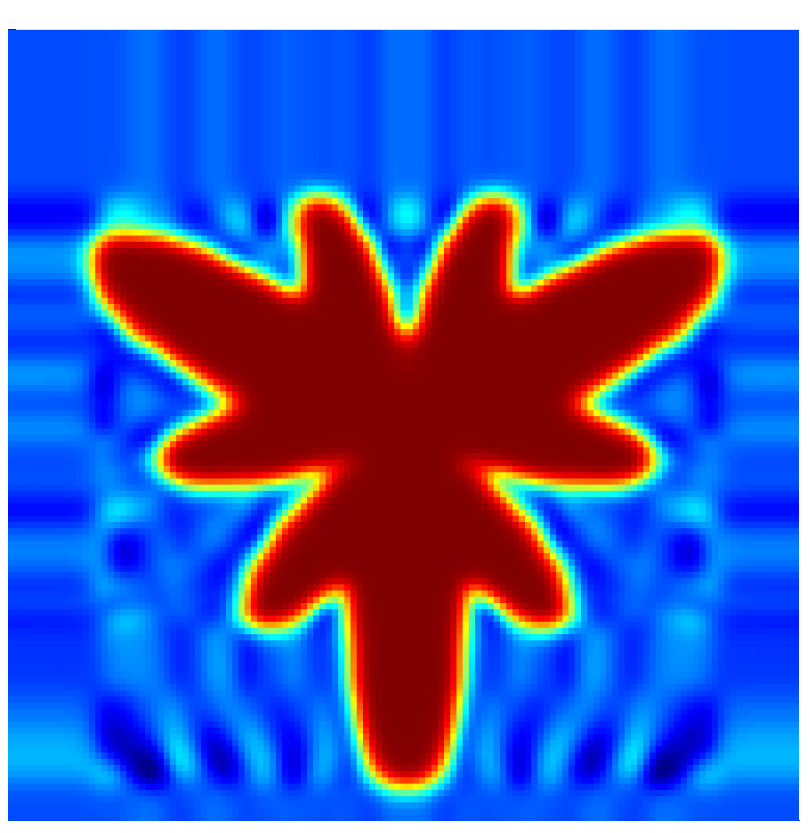} &
\includegraphics[width=1.2in,height=1.2in]{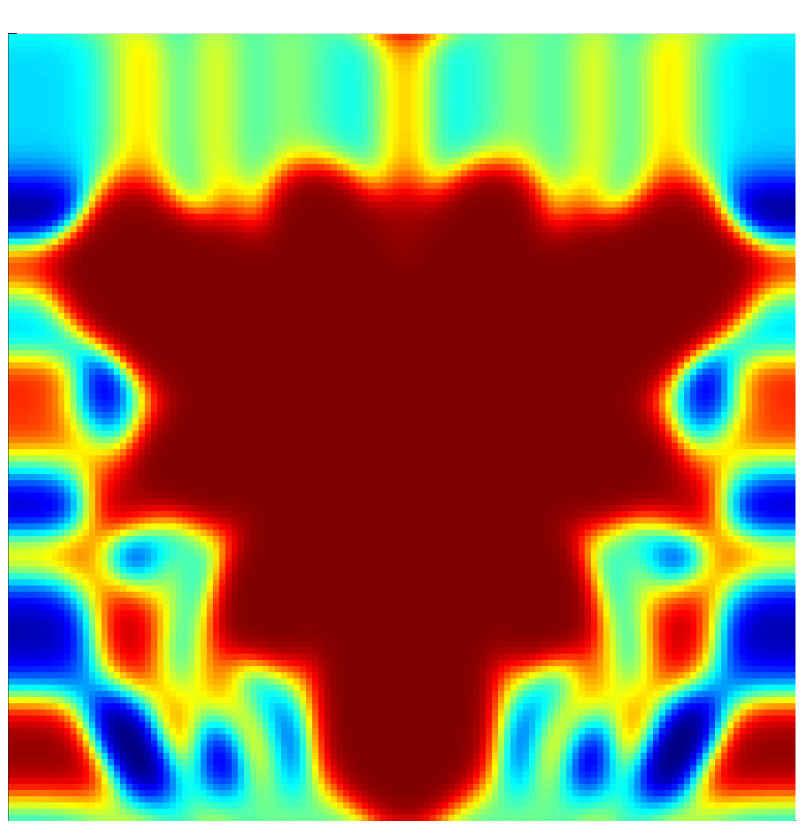} &
\includegraphics[width=1.2in,height=1.2in]{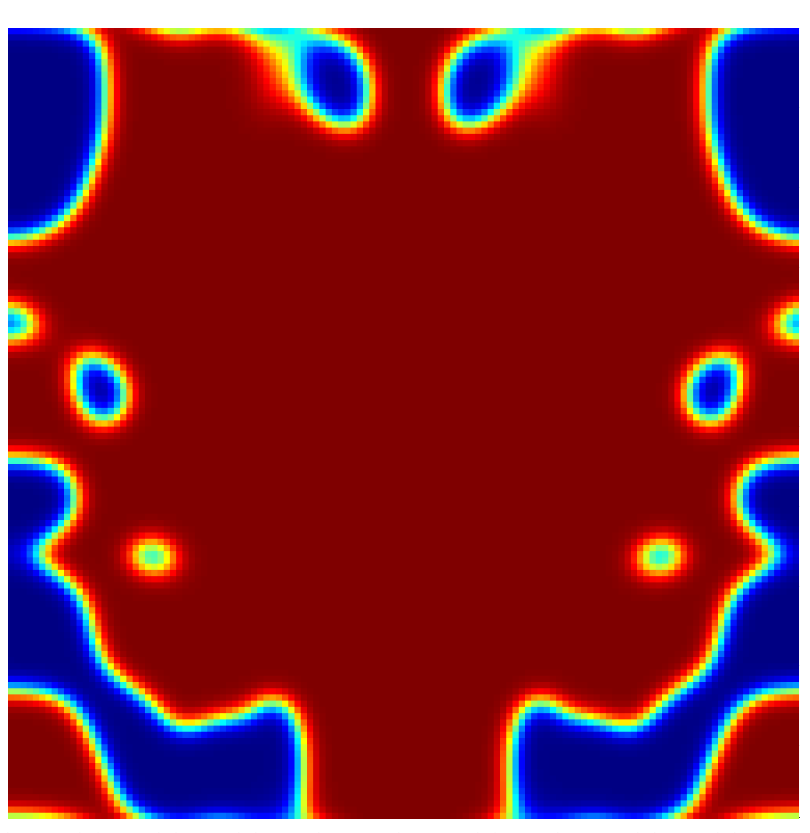}\\
\rotatebox{90}{ADGN-ST} &
\includegraphics[width=1.2in,height=1.2in]{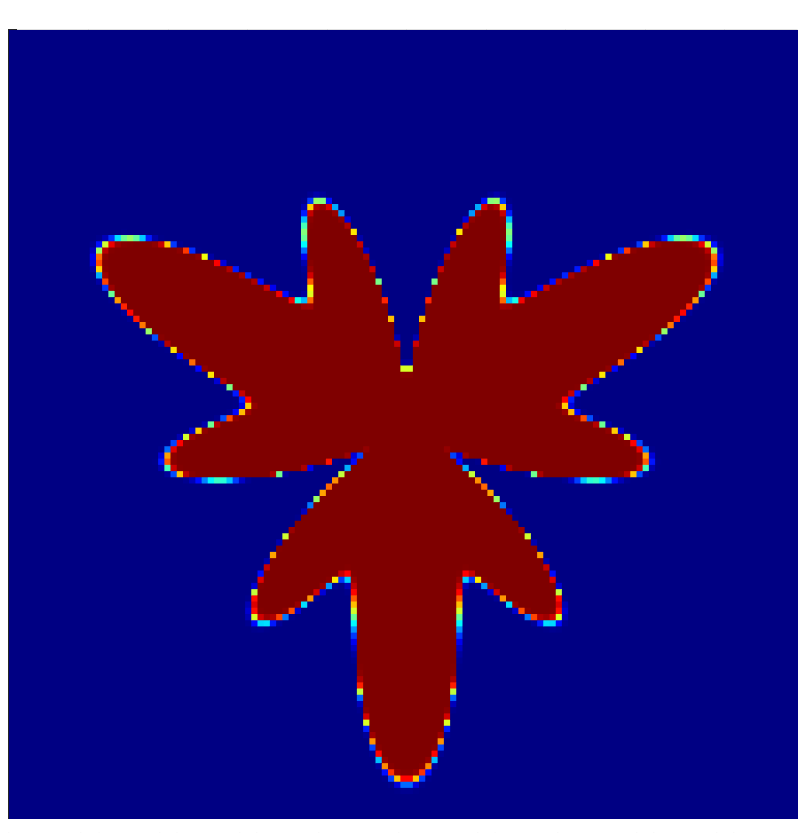} &
\includegraphics[width=1.2in,height=1.2in]{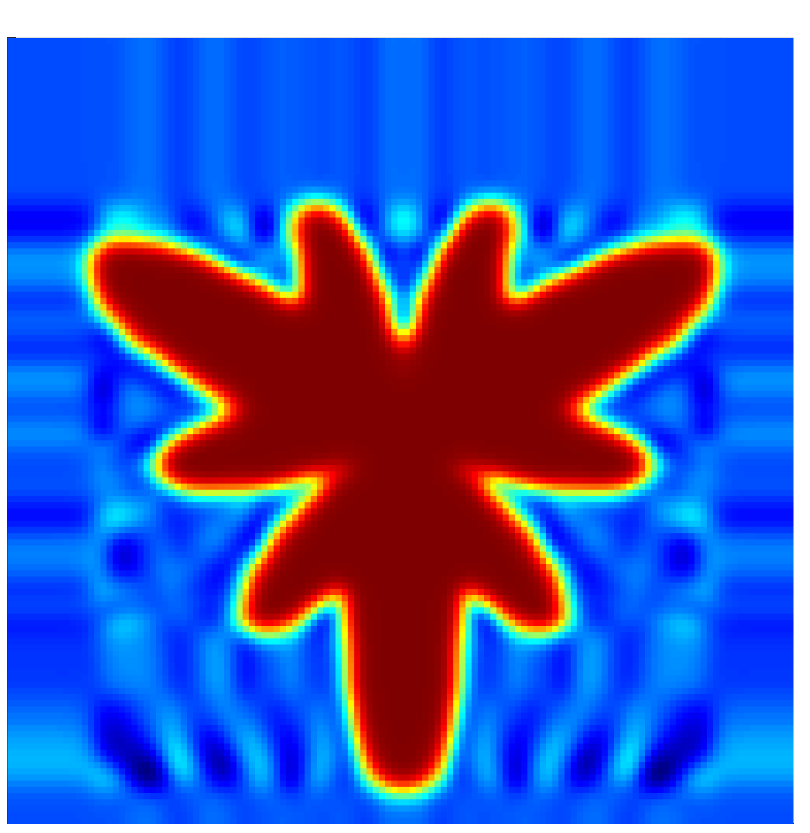} &
\includegraphics[width=1.2in,height=1.2in]{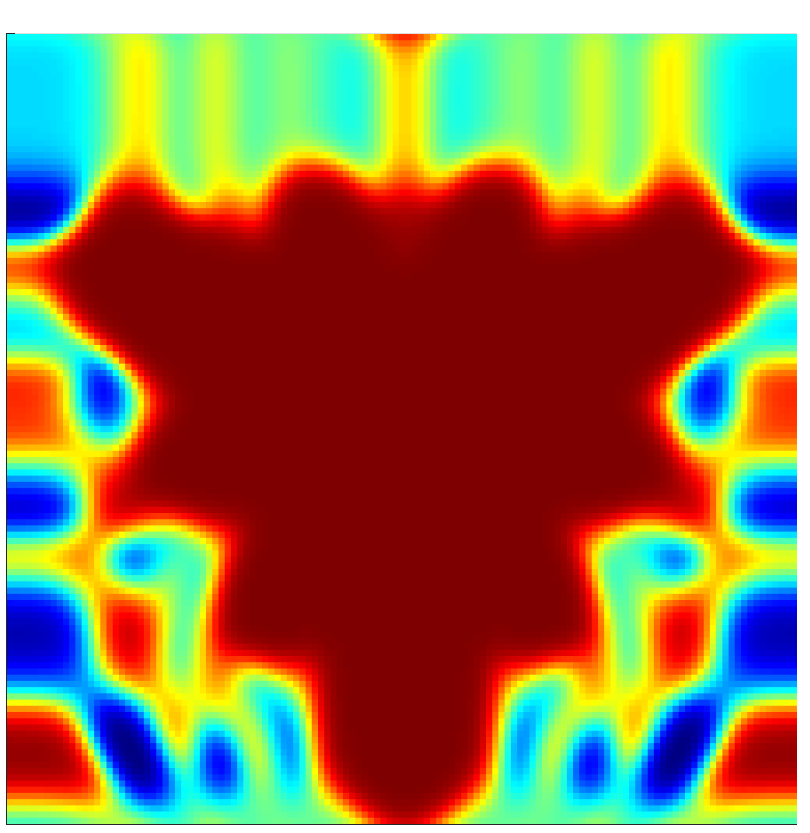} &
\includegraphics[width=1.2in,height=1.2in]{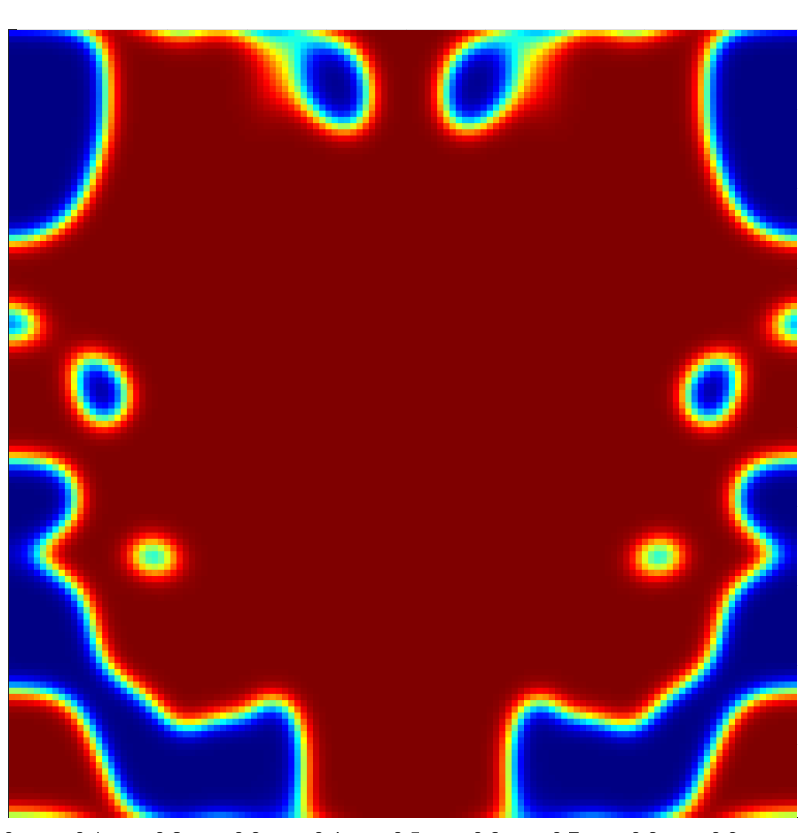} \\
\end{tabular}
\caption{Comparisons of ADRSVD-LT, ADGN-LT, ADRSVD-ST and ADGN-ST: dynamic evolution of butterfly-shaped solutions for Eq.~\eqref{eq4.4} with $\hat{M} = 128$.
Top two rows: Lie-Trotter splitting; Bottom two rows: Strang splitting.}
\label{fig8}
\end{figure}

The temporal evolution of these two initial shapes are shown in Figs.~\ref{fig7} and \ref{fig8}, respectively. 
Fig.~\ref{fig7} illustrates the dynamic evolution of numerical solutions for Eq.~\eqref{eq4.4} with the star-shaped initial condition under spatial grid numbers $\hat{N} = 128$ and the time step size $\tau$ is $0.01$, where the first two rows present the results of ADRSVD-LT and ADGN-LT, and the last two rows show those of ADRSVD-ST and ADGN-ST. Each column corresponds to the phase-field distribution at $T = 0,10,50,100$, respectively. The initially star-shaped interface exhibits a gradual morphological evolution: the tip regions (high-curvature areas) contract inward, while the gap regions (low-curvature areas) expand outward, ultimately approaching a circular morphology-a behavior consistent with curvature-driven flow dynamics. Numerical results demonstrate that both splitting schemes yield solutions with high consistency in interface positioning, curvature distribution, and energy dissipation characteristics.

Fig.~\ref{fig8} presents numerical simulations of the dynamic evolution for Eq.~\eqref{eq4.4} with butterfly-shaped initial condition under spatial grid numbers $\hat{N} = 128$, where the first two rows present the results of ADRSVD-LT and ADGN-LT, and the last two rows show those of ADRSVD-ST and ADGN-ST. Each column correspondsto phase-field distributions at $T = 0, 1, 5, 10$. The initial butterfly-shaped interface evolves with high-curvature regions contracting inward and low-curvature regions expanding outward, leading to a smoothed morphology over time. Results demonstrate consistent accuracy and stability in capturing interface evolution for both splitting schemes, while the adaptive low-rank methods effectively preserve phase-field quality without significant numerical dissipation or interface distortion, validating their robustness in handling complex geometric initial-value problems.

\section{Concluding remarks}
\label{sec5}

This work proposes a splitting-based RDLR approximation 
that unifies computational efficiency and numerical stability. We thoroughly consider the characteristics of \eqref{eq1.1}, dividing the problem into a stiff linear component and a nonstiff nonlinear component. For the stiff part, exponential integrators overcome explicit step size constraints; for the nonlinear part, RDLR approximation enables efficient resolution while bypassing high-dimensional computational bottlenecks. We establish complete integration frameworks for Lie-Trotter and Strang splittings, respectively, with 
extending to the rank-adaptive schemes. Through numerical experiments on the Allen-Cahn equation and the DRE, we verify the convergence orders and computational advantages of the proposed method. Furthermore, we demonstrate its robustness through extensive simulations of high-curvature initial interfaces over extended time periods.
Based on this work, three future research directions are suggested:
\begin{itemize}
\item[(i)] {Extension to tensor manifolds to tackle higher-dimensional stiff PDEs.} 
\item[(ii)]{Establishment of a systematic error analysis framework for the proposed low-rank numerical scheme.}
\item[(iii)] {The current method, utilizing truncated SVD for low-rank approximation, is effective when initial value perturbations are sufficiently small. However, truncated SVD may become unstable or inaccurate when initial perturbations are large. Future work needs to investigate and develop more robust low-rank approximation techniques to effectively handle such more challenging scenarios.}       
\end{itemize}

\section*{Acknowledgments}
\addcontentsline{toc}{section}{Acknowledgments}
\label{sec6}
\textit{This work is supported by the National Natural Science Foundation of China (12401536), Sichuan Science and Technology Program (2024NSFSC0441) and Key Laboratory of Large-Scale Electromagnetic Industrial Software, Ministry of Education (EMCAE202502).}

\section*{Data availability}
The datasets generated and analyzed during the current study are not publicly available but are available from the authors on reasonable request.

\section*{Conflict of interests}
The authors declare that they have no known competing financial interests or
personal relationships that could have appeared to influence the work reported in this paper.

\appendix
\section*{Appendix}
\addcontentsline{toc}{section}{Appendix}
\label{appendix}
\setcounter{algorithm}{0}
\renewcommand{\thealgorithm}{A.\arabic{algorithm}}

For the clarity, this appendix collects the pseudocode of some algorithms in Sections \ref{sec2} and \ref{sec3}.

\begin{algorithm}[ht]
	\caption{Dynamical Rangefinder}
	\label{alg:dyn_rangefinder}
	$\texttt{DynamicalRangefinder}(N_0, \tau, F, r, p, q)$
	\begin{algorithmic}[1]
		\REQUIRE Initial value \( N_0 \in \mathbb{R}^{m \times n} \), step size \( \tau > 0 \), target rank \( r > 0 \), oversampling parameters \( p \geq 0 \), power iterations $ q \geq 0 $
		\STATE Generate Gaussian Sketch Matrix $\Omega\in\mathbb{R}^{n\times(r+p)}$
		\STATE Find \( B(t_0+\tau) \) by solving
		\[
		\dot{B}(t) = F\bigl(t, B(t) W^{\top}\bigr) \Omega, \quad B(t_0) = N_0 \Omega
		\]
		via an explicit fourth-order Runge-Kutta method with 10 substeps.
		\STATE Extract initial orthonormal basis \( Q = \text{orth}(B(t_0+\tau)) \) 
		\FOR {$i = 0$ \textbf{to} $q$}  
		\STATE Compute \( C(t_0+h) \) by solving  
		\[
		\dot{C}(t) = F\bigl(Q C(t)^{\top}\bigr)^{\top} Q, \quad C(t_0) = N_0^{\top} Q
		\]
		via an explicit fourth-order Runge-Kuttamethod with 10 substeps.
		\STATE Extract the orthonormal basis \( Q = \text{orth}(C(t_0+\tau)) \)   
		\STATE  Compute \( B(t_0+\tau) \) by solving   
		\[
		\dot{B}(t) = F\bigl(B(t) Q^{\top}\bigr) Q, \quad B(t_0) = N_0 Q
		\] 
		via an explicit fourth-order Runge-Kuttamethod with 10 substeps.
		\STATE Extract the orthonormal basis \( Q = \text{orth}(B(t_0+\tau)) \)  
		\ENDFOR
		\RETURN \( Q \in \mathbb{R}^{m \times (r+p)} \) 
	\end{algorithmic}
\end{algorithm}

\begin{algorithm}[H]
	\caption{Dynamical Randomized SVD} 
\label{alg:drsvd}
$\texttt{DynamicalRandomizedSVD}(N_0, \tau, F, r, p, q)$
\begin{algorithmic}[1]
	\REQUIRE Initial value $ N_0 \in \mathbb{R}^{m \times n} $ and $N_{0}\approx U_{0}S_{0}V_{0}^{\top}\in\mathcal{M}_{r}$, step size $ \tau > 0 $, target rank $ r > 0 $, oversampling $ p \geq 0 $, power iterations $ q \geq 0 $
	
	\STATE Compute $Q_\tau = \texttt{DynamicalRangefinder}(N_0, \tau, F, r, p, q)$
	
	\STATE Basis augmentation: $Q = \text{orth}([U_0, Q_\tau])$
	
	\STATE Compute $ C(t_0+\tau) $ by solving the reduced $ n\times(2r + p) $ differential equation on subspace
	\begin{align*}
		\dot{C}(t) = F(t,QC(t)^{\top})^{\top} Q, \quad C(t_0) = N_0^{\top} Q
	\end{align*}
	via an explicit fourth-order Runge-Kuttamethod with 10 substeps.
	
	\STATE Compute a SVD of $ C(t_0+\tau)^{\top} = \tilde{U}_{\tau} \Sigma_\tau V_\tau^{\top} $ and $N(t_{0}+\tau)=Q\tilde{U}_{\tau}\Sigma_{\tau}V_{\tau}^{\top}$
	\STATE Truncate the SVD to the target rank: $U_{1}\Sigma_{1}V_{1}^{\top}=\mathcal{T}_{r}(Q\widetilde{U}_{\tau}\Sigma_{\tau}V_{\tau}^{\top})$
	\RETURN {$N_1 = U_1 \Sigma_1 V_1^{\top}$}
\end{algorithmic}
\end{algorithm}

\begin{algorithm}[p]
\caption{Dynamical Generalised Nyström}
\label{alg:dgn}
$\texttt{DynamicalGeneralisedNyström}(N_0, \tau, F, r, p, \ell, q)$
\begin{algorithmic}[1]
	\setlength{\baselineskip}{0.85\baselineskip}  
	\REQUIRE Initial value \( N_0 \in \mathbb{R}^{m \times n} \) and \(N_{0}\approx U_{0}S_{0}V_{0}^{\top}\in\mathcal{M}_{r}\), step size \( \tau > 0 \), target rank \( r > 0 \), oversampling parameters \( p \geq 0 \) and \( \ell \geq 0 \), power iterations \( q \geq 0 \)
	
	\STATE Compute in parallel: 
	\[
	\tilde{Q}_1 = \texttt{DynamicalRandomizedRangefinder}(N_0, \tau, F, r, p, q)
	\]
	\[
	\tilde{Q}_2 = \texttt{DynamicalRandomizedCo-rangefinder}(N_0, \tau, F, r, p+\ell, q)
	\]
	\STATE Augment the two bases: 
	\[
	Q_1 = \text{orth}([U_0, \tilde{Q}_1]), \quad Q_2 = \text{orth}([V_0, \tilde{Q}_2])
	\]
	
	\STATE Compute \( B(\tau) \) by solving the \( m \times (r + p + \ell) \) differential equation (row space evolution):
	\[
	\dot{B}(t) = F(t,B(t)Q_2^{\top})Q_2, \quad B(t_0) = A_0 Q_2
	\]
	via an explicit fourth-order Runge-Kuttamethod with 10 substeps.
	\STATE Compute \( C(\tau) \) by solving the \( n \times (r + p) \) differential equation (column space evolution):
	\[
	\dot{C}(t) = F(t,Q_1C(t)^{\top})^{\top} Q_1, \quad C(t_0) = A_0^{\top} Q_1
	\]
	via an explicit fourth-order Runge-Kuttamethod with 10 substeps.
	\STATE Compute \( D(\tau) \) by solving the \( (r + p) \times (r + p + \ell) \) differential equation (interaction coupling terms):
	\[
	\dot{D}(t) = Q_1^{\top} F(t,Q_1D(t)Q_2^{\top})Q_2, \quad D(t_0) = Q_1^{\top} A_0 Q_2
	\]
	via an explicit fourth-order Runge-Kuttamethod with 10 substeps.
	\STATE Conduct a SVD of \( D(t_0+\tau) = U_\tau S_\tau V_\tau^{\top} \)
	\STATE Truncate to the target rank: \( \mathcal{T}_r(D(t_0+\tau)) = U_rS_rV_r^{\top} \) 
	\STATE Compute QR decompositions and assemble:
	\[
	UR_{\mathrm{l}}=B(t_{0}+\tau)U_{r}, \quad VR_{2}=C(t_{0}+\tau)^{\top}V_{r} \quad \leadsto \quad S = R_1 S_r R_2^{\top}
	\]
	\RETURN {\( N_1 = U S V^{\top} \)}
\end{algorithmic}
\end{algorithm}

\begin{algorithm}[p]
\caption{Adaptive Dynamical Rangefinder}
\label{alg:adaptive_rangefinder}
$\texttt{AdaptiveDynamicalRangefinder}(N_0, \tau, F, \varphi, \beta)$
\begin{algorithmic}[1]
	\REQUIRE Initial value \( N_0 \in \mathbb{R}^{m \times n} \), step size \( \tau > 0 \), tolerance \( \varphi  > 0 \), max failure probability \( \beta > 0 \)
	
	\STATE Set \( \kappa = -\left\lfloor \frac{\log(\beta)}{\log(10)} \right\rfloor \) and \( \varepsilon = \sqrt{\frac{\pi}{2}} \cdot \frac{\varphi }{10} \) 
	\STATE Initialize a Gaussian random matrix \( \Omega \in \mathbb{R}^{n \times \kappa } \)
	\STATE Compute \( B(t_0+\tau) \) by solving the \( m \times \kappa \) differential equation:
	\[
	\dot{B}(t) = F(t,B(t)\Omega^{\top} \Omega)^{-1} \Omega^{\top}) \Omega, \quad B(t_0) = N_0 \Omega
	\]
	\STATE Extract the basis: \( Q_\tau = \text{orth}(B(\tau)) \) 
	via an explicit fourth-order Runge-Kuttamethod with 10 substeps.
	\STATE Set \( E = \infty \) and \( j = 1 \)
	
	\WHILE{\( E > \varepsilon \)} 
	\STATE Generate a new Gaussian random matrix \( \Omega \in \mathbb{R}^{n \times \kappa} \)
	\STATE Compute \( B(t_0+\tau) \) by solving the \( m \times \kappa \) differential equation:
	\[
	\dot{B}(t) = F(t,B(t)\Omega^{\top} \Omega)^{-1} \Omega^{\top}) \Omega, \quad B(t_0) = N_0 \Omega
	\]
	via an explicit fourth-order Runge-Kuttamethod with 10 substeps.
	\STATE Compute residual: \(\tilde{B}=B(t_0+\tau)-Q_{\tau}Q_{\tau}^{\top}B(t_0+\tau)\)
	\STATE Update \(E=\max_{i=1,...,\kappa}\left\|\tilde{B}[:,i]\right\|\)
	\STATE Set \( j = j + 1 \) and update the basis: \( Q_\tau = \text{orth}([Q_\tau, \tilde{B}]) \)
	\ENDWHILE
	\RETURN {\( Q_\tau \in \mathbb{R}^{m \times j\kappa} \)}
\end{algorithmic}
\end{algorithm}

%
%
%
%




%
%
%
\bibliographystyle{elsarticle-num}   
\bibliography{reference}             

\end{document}